\documentstyle{article}

\font\tenmsb=msbm10
\font\sevenmsb=msbm7
\font\fivemsb=msbm5

\newfam\msbfam
\textfont\msbfam=\tenmsb
\scriptfont\msbfam=\sevenmsb
\scriptscriptfont\msbfam=\fivemsb
\def\Bbb#1{{\fam\msbfam #1}}

\font\teneufm=eufm10
\font\seveneufm=eufm7
\font\fiveeufm=eufm5
\newfam\eufmfam
\textfont\eufmfam=\teneufm
\scriptfont\eufmfam=\seveneufm
\scriptscriptfont\eufmfam=\fiveeufm
\def\frak#1{{\fam\eufmfam\relax#1}}

\newcommand\qed{{\hspace*{\fill}Q.E.D.\vskip12pt plus 1pt}}
\newcommand\sD{{\cal D}}
\newcommand\sE{{\cal E}}
\newcommand\sF{{\cal F}}
\newcommand\sG{{\cal G}}
\newcommand\sH{{\cal H}}
\newcommand\sJ{{\cal J}}
\newcommand\sK{{\cal K}}
\newcommand\sL{{\cal L}}
\newcommand\sM{{\cal M}}
\newcommand\sO{{\cal O}}
\newcommand\sP{{\cal P}}
\newcommand\sQ{{\cal Q}}

\newcommand\gra{\alpha}
\newcommand\grb{\beta}

\newcommand\grk{\kappa}

\newcommand\vphi{\varphi}

\newcommand\grs{\sigma}

\newcommand\rat{{\Bbb Q}}

\newcommand\comp{{\Bbb C}}
\newcommand\zed{{\Bbb Z}}

\newcommand\pn[1]{{\Bbb P}^{#1}}
\newcommand\proj[1]{{\Bbb P}({#1})}
\newcommand\pnsheaf[2]{{{\cal O}_{\pn{#1}}({#2})}}
\newcommand\pnpair[2]{{(\pn {#1},\pnsheaf{#1}{#2})}}

\newcommand\oline[1]{{\overline {#1}}}

\def\rank{{\mathop{\rm rank}\nolimits}}
\def\cod{{\mathop{\rm cod}\nolimits}}
\def\min{{\mathop{\rm min}\nolimits}}
\def\dim{{\mathop{\rm dim}\nolimits}}

\newcommand\Pic[1]{{{\rm Pic(}{#1}{\rm )}}}

\newcommand\proof{{\noindent\bf Proof.\ }}

\newtheorem{theorem}{Theorem}[section]
\newtheorem{lemma}[theorem]{Lemma}
\newtheorem{corollary}[theorem]{Corollary}

\newtheorem{prop}[theorem]{Proposition}
\newtheorem{definition}[theorem]{Definition}
\newtheorem{re}[theorem]{Remark}
\newtheorem{pargrph}[theorem]{}
\newtheorem{examp}[theorem]{Example}

\newenvironment{rem*}{\begin{re}\em}{\end{re}}
\newenvironment{example}{\begin{examp}\em}{\end{examp}}
\newenvironment{def*}{\begin{definition}\em}{\end{definition}}

\newenvironment{prgrph*}[1]{\indent\begin{pargrph}{\bf #1.}\em\
}{\end{pargrph}}

\begin{document}

\title{Projections from Subvarieties}
\author{Mauro C. Beltrametti, Alan Howard, Michael Schneider and
         Andrew J. Sommese}

\date{April 2, 1998}
\maketitle

\bigskip

\tableofcontents

\section{Introduction}Let $X\subset \pn N$ be an $n$-dimensional
 connected projective submanifold of projective space.  Let
$p : \pn N\to \pn {N-q-1}$ denote the projection from a linear
$\pn q\subset \pn N$.
Assuming that $X\not\subset \pn q$ we have the induced rational
mapping $\psi:=p_X: X\to \pn {N-q-1}$. This article
started as an attempt to understand the structure of this mapping when
$\psi$ has a lower dimensional image.  In this case of necessity
we have $Y := X\cap \pn q$ is nonempty.

The special case when $Y$ is a point is very classical: $X$ is a linear
subspace of $\pn N$.  The case when $q=1$ and $Y=\pn q=\pn 1$ was settled
for surfaces by the fourth author \cite{So1} and by Ilic \cite{Ilic} in
general.
Beyond this even the special case when $q\ge 2$ and $Y=\pn q$ is open.

 We have found it convenient to
study a closely related question, which includes many special cases
including the
case when the center of the projection $\pn q$ is contained in $X$.

\noindent{\bf Problem.} Let $Y$ be a proper connected $k$-dimensional
projective submanifold of an $n$-dimensional projective manifold $X$. 
Assume that $k>0$. Let $L$
be a very ample line bundle on $X$  such that
$ L\otimes \sJ_Y$ is spanned by global sections, where $\sJ_Y$ denotes 
the ideal sheaf of $Y$ in
$X$.  Describe the structure of
$(X,Y,L)$ under the additional assumption that the image of $X$ under the
mapping $\psi$ associated to $| L\otimes \sJ_Y|$ is lower dimensional.

 Let us describe our progress on this problem.

In \S \ref{LUBSec} we study upper and lower bounds for the dimensions
of the spaces of sections of powers $tL$ of a very ample line bundle $L$ on a
projective manifold $X$.  The need for such bounds arises naturally when
we consider line bundles  which are multiples of a very ample line bundle.
One general result Proposition (\ref{lowerbound}) gives  
an upper bound for an integer $t_0$ such that for $t\ge t_0$,
$h^0(tL \otimes \sJ_Y) > 0$.

In \S \ref{StructureSec} we prove a number of general results.  For
example, Theorem
(\ref{Dimension})
shows that $\dim \psi(X)\ge n-k-1$ with equality only if $Y$ is a complete
intersection in $X$.  In particular   Corollary (\ref{Pk}), shows that
 if $Y$ is a linear $\pn k$, then
$\dim \psi(X)\ge n-k-1$ with equality only if $(X,L)\cong \pnpair n 1$.
Proposition (\ref{BetterB}) further shows that if $\rank\Pic X=1$ and $\dim
\psi(X)\ge n-k$
then
$\displaystyle \dim \psi(X)\geq n-k+\frac{k}{n-k}-1$.  Theorem (\ref{ThmB})
shows that if
$Y$ is a $\pn k$ (or more generally a projective manifold whose algebraic
cohomology
is the same as $\pn k$ up to dimension $2(n-k)$), then if
$\dim \psi(X)\ge n-k$, it follows that $\dim \psi(X)\ge k$.  In particular
except for
known examples, we have for a wide range of $Y$ including $\pn k$, that
$\ \displaystyle \dim \psi(X)\ge  \frac{\dim X}{2}$.

In \S \ref{ExSec} we give a number of examples showing that the dimensions
allowed by the examples
do occur.  Of particular interest is Example (\ref{Man7}).  This example
consists for each positive integer $n$
of an infinite sequence of projective $n$-folds in $\pn {2n-1}$ which
contain a linear $\pn{n-1}$. All
degrees of $X$ that are allowed by theory occur.

In \S \ref{Divisorial} we specialize to the case when $Y$ is a divisor. We
study bundles of the form $tL-Y$ where $t$ is near $\delta:=\deg Y$. One 
result, Theorem (\ref{Properties}), implies that if $\delta > 1$ then 
$|\delta L-Y|$ 
gives a birational map, which is in fact very ample if 
$2n\ge  \dim \Gamma(L) +1$.

In \S \ref{Linear} we restrict to the case when $Y$ is
a linear $\pn k$ and show, among other things, that $\dim Z \geq n-k$
except when $X$ is a hypersurface in $\pn{n+1}$. In \S \ref{Linear1}  
we restrict further to the special case when $Y$ is a linear $\pn{n-1}$.
In this case $\psi$ is a morphism.  Remmert-Stein factorize $\psi = s\circ
\phi$ with $\phi :X\to Z$
a morphism with connected fibers onto a normal projective variety $Z$, and
with $s$ a finite morphism.
We know that except for known examples, if $\dim\psi(X)<\dim X$  then
$\dim\psi(X)=n-1$.
We show that $Z$ is very well behaved (Cohen-Macaulay, $\rat$-factorial,
$\Pic Z\cong \zed$).  Moreover we examine the possible degrees of $s$ and
use adjunction theory to
classify the possible $(X,L)$ for extreme values of this degree.

We would like to thank Frank-Olaf
Schreyer for his very helpful explanation of how Castelnuovo theory 
gives lower bounds
for the dimensions of spaces of sections of powers of very ample line 
bundles. 

The research in this article was carried out in Bayreuth,
the University of Notre Dame,
and two sessions of the RiP program at Oberwolfach.  All the authors are
indebted to the Volkswagen Stiftung, whose generosity allowed us to
work together in such an ideal setting.  The fourth author thanks the
Alexander  von Humboldt Stiftung  for their generous support.  

The final stages of this article were developed during a three-week stay at
Oberwolfach
in  the summer of 1997.   Within a few weeks after we separated,
Michael Schneider died in a climbing accident.  The three remaining
authors dedicate this work to his memory. He was our friend and
colleague, a person of vibrant energy, keen intelligence, and
generous spirit.  We feel a deep sense of loss, but are grateful to
have had our lives and work enriched by his presence.

\section{Background material}\label{Background}
\addtocounter{subsection}{1}\setcounter{theorem}{0}
We work over the complex numbers $\comp$. Through the paper we deal with
projective varieties $V$. We denote by $\sO_V$ the structure sheaf of $V$
and by $K_V$ the canonical bundle, for $V$ smooth.
For any coherent sheaf $\sF$ on $V$, $h^i(\sF)$ denotes the complex
dimension of
$H^i(V,\sF)$.

Let $\sL$ be a line bundle on $V$. The line bundle $\sL$ is said to be {\em
numerically effective {\rm
(}nef}, for short) if $\sL \cdot C \geq 0$ for all effective curves $C$ on
$V$. $\sL$ is said to be
{\em big} if $\grk(\sL) = \dim V$, where $\kappa(\sL)$ denotes the Kodaira
dimension of $\sL$. If $\sL$ is
nef then this is equivalent to $c_1(\sL)^n > 0$, where $c_1(\sL)$ is the
first Chern class of $\sL$
and $n= \dim V$.

\begin{prgrph*}{Notation}The notation used in this paper is standard from
algebraic geometry. In particular, $\approx$ denotes
 linear equivalence of line bundles. For a line bundle $\sL$ on a compact
complex space $V$,
  $\chi(\sL) := \sum_{i}(-1)^i h^i(\sL)$ denotes  the Euler characteristic, 
and $|\sL|$ denotes the complete linear system associated with a
line bundle. We say that $\sL$ is spanned if it is spanned  at all points
of $V$ by $\Gamma(\sL)$.

For a compact connected projective manifold $V$,
$h^{2j}(V,\rat)_{\rm alg}$ denotes the  dimension of the vector
subspace   $H^{2j}(V,\rat)_{\rm alg}$ of $H^{2j}(V,\rat)$ dual under
Kronecker duality to the vector subspace of $H_{2j}(V,\rat)$
spanned by the $j$-dimensional  algebraic subvarieties of $V$.

We denote  the ideal
sheaf of an irreducible subvariety $A$ of a variety $V$ by
$\sJ_{A/V}$ (or simply $\sJ_A$ when no confusion can result). For smooth
$A$ contained
in the smooth locus of $V$, $N_{A/V}$ denotes the normal
bundle of $A$ in $V$.

Line bundles and divisors are used with little (or no)
distinction. Hence we shall freely switch between the multiplicative 
and the additive notation.
\end{prgrph*}

\begin{prgrph*}{Conductor formula} \label{ConductorF} Let $V$ be a
connected projective manifold of dimension $n$. Let $L$ be a very
ample  line bundle on $V$ of degree $d:=L^n$ with $|L|$ embedding
$V$ into $\pn N$. Then the classical conductor formula  for the
canonical bundle states that
$$\Delta\in |(d-n-2)L-K_V|,$$
where $\Delta$ is the double point divisor of a projection of $V$ from
$\pn N$ to  $\pn{n+1}$ (in the degenerate cases when
$N=n$ or $n+1$, $\Delta$ is taken to be the empty divisor).
In particular the line bundle $(d-n-2)L-K_V$ is spanned since given any
point of $V$
a generic projection can be chosen
  with the point not in the double point divisor of the projection
 (see \cite{Zar0} and \cite[p. 71]{Zar1}).
\end{prgrph*}

The following standard lemma is basic (see also \cite[(3.1.8)]{Book}).

\begin{lemma}\label{BasicLemma}Let $V$ be an irreducible
 normal projective variety with $\Pic V\cong \zed$. Let $g:V\to Z$ be a
surjective
morphism of $V$ to a projective variety $Z$.
   Either
$g$ is  a finite morphism or $g(V)$ is a point. The same conclusion
holds for
$V\cong \pn n$ and any holomorphic map to a compact complex space $Z$.
\end{lemma}
\proof Assume that $g$ is not finite and doesn't map $V$ to a point. If
$Z$ is projective, then the pullback of an ample line bundle cannot be 
ample and thus
we see that $\Pic V\not \cong \zed$.  Thus we can assume that $V\cong \pn
n$ and
$Z$ is not necessarily projective.

Note that $\dim g(\pn n)=n$.  If not let $F$ denote a general fiber.
Since
it is smooth it would have  trivial normal bundle.  This contradicts the
ampleness of the
tangent bundle of $\pn n$.

Let $F$ denote a positive dimensional fiber. There is a complex
neighborhood $U$ of $F$
which maps generically one-to-one to a Stein space.  Since $\pn n$ is
homogeneous we
have that the translates of $F$ fill out an open set.
Since the map
$g_U$ must map these positive  dimensional subspaces to points we have the
contradiction
that $\dim g(U)=\dim g(\pn n)<n$.
\qed

We also need the following general fact.

\begin{lemma}\label{CM} Let $f:X\to Y$ be a surjective proper map between
normal varieties.
Assume that $X$ is Cohen-Macaulay and all fibers of $f$ are equal
dimensional. Then $Y$ is
Cohen-Macaulay.
\end{lemma}
\proof Note that a Cartier divisor on a Cohen-Macaulay variety is
Cohen-Macaulay and that if we slice with $\dim X-\dim Y$ sufficiently
ample divisors, then the restriction of the map to the slice is finite by a
well known theorem of Hironaka \cite[(2.1)]{Hironaka}.  Since a general
hyperplane section of a normal variety is normal by Seidenberg's
theorem,  we can assume without loss of generality
that $f$ is finite. Let $n:=\dim X=\dim Y$. By using
\cite[II, (7.6)]{Hartshorne} we are reduced to showing that for any
locally  free coherent sheaf $\sE$ on $Y$, we have
\begin{equation}\label{vanish}
h^i(\sE(-q))=0\;\;{\rm for}\; i<n\;{\rm and}\; q\gg 0,
\end{equation}
where $\sF(t)$ for a coherent sheaf $\sF$ means $\sF\otimes H^{\otimes t}$
for a fixed ample
line bundle $H$ on $Y$.

Since $f$ is finite the pullback of an ample line bundle is ample.
Given  a coherent sheaf $\sG$ on $X$, $\sG(t)$ means $\sG$ twisted by
the $t$-th power of the pullback of $H$.
Hence $(f^*\sE)(t)=f^*(\sE(t))$.

Now since $X$ is Cohen-Macaulay we have
$$h^i((f^*\sE)(-q))=0\;\;{\rm for}\;i<n\;{\rm and}\;q\gg 0.$$
By the Leray spectral sequence, the projection formula and vanishing of
higher direct images
we obtain
\begin{equation}\label{iso}
h^i((f^*\sE)(-q))=h^i(\sE(-q)\otimes f_*\sO_X).
\end{equation}
Since both $X$ and $Y$ are normal we can use the trace mapping from
 $f_*\sO_X\to \sO_Y$ to see that the
exact sequence
$\ 0\to \sO_Y\to f_*\sO_X\to \sM\to 0\ $
splits, where $\sM$ denotes the quotient bundle.  Thus $f_*\sO_X\cong\sO_Y\oplus \sM$.
Thus by combining
(\ref{vanish}) and (\ref{iso}) we get
$$0=h^i((f^*\sE)(-q))=h^i(\sE(-q))+h^i(\sE(-q)\otimes \sM),$$
for $i<n$ and $q\gg 0$. Then (\ref{vanish}) follows and we are done.\qed

If $X$ is smooth and $f$ is finite we can say more.

\begin{lemma}\label{Piclemma}
Let $f: X \to Y$ be a finite surjective map between projective varieties,
where $X$ is smooth and $Y$ is normal.
Then $Y$ is Cohen-Macaulay and $(\deg f)$-factorial.  Moreover, if $-K_X$
is nef
and big, the induced map of
${\rm Pic}(Y) \to {\rm Pic}(X)$ is injective.
\end{lemma}
\proof
The fact that $Y$ is Cohen-Macaulay was proved in the previous lemma.  To
see that it is $(\deg f)$-factorial, let $D$ be a Weil
divisor on $Y$.  Since $X$ is smooth, $f^*D$ is a Cartier divisor.  We
construct a Cartier divisor, ${\rm Norm}(f^*D)$, on $Y$ as
follows: in a small neighborhood $U$ of any smooth point $y$ in $Y$ over
which $f$ is unramified, we define a rational function
by multiplying the functions defining $f^*D$ on the connected components of
$f^{-1}(U)$; and we construct the divisor
determined locally by this construction, first over all smooth points of
$Y$ over which $f$ is unramified, and then (since $Y$
is normal) to all of $Y$ by Riemann extension.

 From the way ${\rm Norm}(f^*D)$ was constructed it is obvious that ${\rm
Norm}(f^*D) = (\deg f) D$.  This shows that $Y$ is
$(\deg f)$-factorial.  In addition, the same construction shows that if
$D$ is a Cartier divisor on $Y$ for which $f^*D$ is trivial, then $(\deg
f)\cdot
D$ is trivial.  In particular, the kernel of the
induced map of ${\rm Pic}(Y) \to {\rm Pic}(X)$ consists entirely of torsion
elements.

Now suppose that $-K_X$ is nef and big.  Then $h^i(\sO_X) = 0$ for $i > 0$.
Using the direct sum decomposition
 $f_*\sO_X\cong\sO_Y\oplus \sM$ from the previous lemma together with the
Leray spectral sequence applied to the finite
map $f$, we see that $h^i(\sO_Y) = 0$ for $i > 0$.  Therefore, ${\rm
Pic}(Y) \cong H^2(Y, \zed)$, and it follows that
${\rm Pic}(Y) \to {\rm Pic}(X)$ is injective unless there is torsion in
$H^2(Y, \zed)$.  We will show this can not occur.

By the universal coefficient theorem, torsion in $H^2(Y, \zed)$ is
equivalent to torsion in $H_1(Y, \zed)$, which in turn
implies the existence of a finite unbranched covering $Y' \to Y$.  Lifting
this to $X$ gives the commutative diagram
$$
\begin{array}{cccl}
X'&\rightarrow&X \\
\downarrow& &\downarrow&f\\
Y'&\rightarrow&Y
\end{array}
$$
where $Y'$ is connected, the vertical arrows are branched coverings, and
the horizontal arrows are unbranched coverings.
Let $m$ be the common sheet number of both of the latter.
It is easy to see that $X'$  consists of a finite number of disjoint
connected components, each mapping isomorphically onto
$X$; for $h^i(\sO_{X'})= 0$ (because $-K_{X'}$ is big and nef), and
$\chi(\sO_{X'}) = m\chi(\sO_X) = m$, where $m$ is the
sheet number. If $A$ is any connected component of $X'$, we thus get a
finite surjective map  $X \cong A \to Y'$.
Arguing as before, we see that $h^i(Y') =0$ for $i >0$, so that
$\chi(\sO_{Y'}) = 1$.  On the other hand, we have
$\chi(\sO_{Y'}) = m\chi(\sO_Y) = m$.
\qed

The following general lemma is well known and follows from the results
in the introduction of \cite{Mumford} (see also \cite[(6.6.1)]{Book}).

\begin{lemma}\label{Mumford} Let $L$ be a very ample line bundle on an
irreducible projective variety, $X$. Let $Y\subset X$ be an irreducible
subvariety of degree $\delta$ relative to $L$, i.e., $\delta=L^{\dim Y}\cdot
Y$. If either $Y$ is smooth or $Y\subset{\rm reg}(X)$ and ${\rm
cod}_XY=1$ then $\sJ_Y(\delta)$ is spanned by global sections,
where $\sJ_Y(\delta)$ denotes the ideal sheaf $\sJ_Y$ of $Y$ in $X$
tensored with
$\delta L$.
\end{lemma}

The following result we need is a ``folklore'' result, for which we
don't know references.

\begin{prop}\label{HartRev} Assume that Hartshorne's conjecture
 {\rm \cite{HartshCI}}, that any connected nondegenerate $n$-dimensional
smooth submanifold $X\subset \pn m$ is  a
complete intersection if $n>\frac{2}{3}m$, is true.
Then each vector bundle $\sE$ on $\pn m$ of rank $r<\frac{m}{3}$
splits into a direct sum of line bundles.\end{prop}
\proof We use induction over $r$. If $r=1$ the assertion is true. So,
let us assume the assertion true for $r-1$.

Since the assertion is independent of twisting, we may assume that $\sE$
is generated by global sections. Take a general section $s\in H^0(\sE)$
and let $X:=V(s)$,  the zero locus of $s$. Then $X$ is smooth and ${\rm
cod}_{\pn m}X=r$. The assumption $r<\frac{m}{3}$ is equivalent to $\dim
X=m-r>\frac{2}{3}m$ and therefore $X$ is a complete intersection in $\pn
m$ by Hartshorne's conjecture. Thus
$$\sE_X\cong N_{X/\pn m}\cong\oplus_{i=1}^r\sO_X(a_i),$$
where $\sE_X$ denotes the restriction of $\sE$ to $X$. Since $N_{X/\pn
m}$ is ample, the $a_i$'s are positive integers and we may assume
$a_1\geq a_2\geq\cdots\geq a_r>0$.

\noindent{\bf Claim.\ } $\sE(-a_1)$ has a section without zeros.

Assuming the Claim true, we get an exact sequence (given by that section)
$$\ 0\to\sO_{\pn m}\to\sE(-a_1)\to\sF\to 0\ ,$$
where the quotient $\sF$ is a rank $r-1$ vector bundle. Then by
induction $\sF$ splits. Therefore $\sE(-a_1)$, and hence $\sE$, splits
into a direct sum of line bundles.

Thus it remains to show the Claim. Note that
$$\sE_X(-a_1)\cong\sO_X\oplus\sO_X(a_2-a_1)\oplus\cdots\oplus
\sO_X(a_r-a_1),$$ where $a_i-a_1\leq 0$ for $2\le i\le r$. Let $\grs\in
H^0(\sO_X)$ be a section of $\sE_X(-a_1)$ with no zeros. The
obstruction to extending $\grs$ to a formal neighborhood $\widehat{X}$ of $X$
belongs to $H^1(X,S^t(N^*)\otimes\sE_X(-a_1))$, where $N:=N_{X/\pn m}$.
Since $S^t(N^*)\otimes\sE_X(-a_1)$ is a direct sum of negative line
bundles, we have $H^1(S^t(N^*)\otimes\sE_X(-a_1))=0$, $t\geq 1$, by Kodaira
vanishing.
Thus we conclude that there exists a section $\widehat{\grs}\in
H^0(\widehat{X},\sE_{\widehat{X}}(-a_1))$ whose restriction to $X$ coincides
with $\grs$. As soon as $\dim X\geq2$ (which is the case since $m\geq 3$),
it is a fairly standard fact, by using
results of Barth \cite[Proposition 4]{Barth} and Griffiths
\cite[Theorems I, III, p. 378, 379]{Griffiths} (see also \cite[p. 226,
227]{HartshorneA}), that $\widehat{\grs}$ extends to a section
$\tau\in H^0(\sE(-a_1))$. Then the restriction $\tau_X$ has no zeros on
$X$. We want to show that $\tau$ has no zeros on $\pn m$. Let
$Y:=V(\tau)$ be the zero locus of $\tau$. If $Y\not= \emptyset$, then
$\dim Y\geq m-{\rm rank}\sE=m-r$. Since $r<\frac{m}{3}$,
we have that $\ \displaystyle \dim(X\cap Y)\geq \dim X+\dim Y-m\geq  m-2r>0$.
Therefore $X\cap Y\not= \emptyset$ in $\pn m$. This contradicts the fact
that the restriction $\tau_X$ has no zeros on $X$.\qed

\section{Lower and upper bounds for $h^0(tL)$}\label{LUBSec}
\addtocounter{subsection}{1}\setcounter{theorem}{0}
We first state some general lower and upper bound formulas for the
number of sections of multiples of a given line bundle $L$.

\begin{lemma}\label{upper} Let $L$ be a big and spanned line bundle on
an irreducible  $n$-dimensional projective variety $X$. Then for
$t\geq 0$ with $d:={\rm deg}_L(X)=L^n$ we have
$$h^0(tL)\leq {{t+n-1}\choose
{n-1}}\frac{dt+n}{n}
=\frac{td+n}{t+n}{{t+n}\choose {n}}$$
\end{lemma}
\proof If $X$ is a curve then clearly the result is true, i.e.,
$\ \displaystyle h^0(tL)\le {t+1-1\choose 1-1}\frac{dt+1}{1}=dt+1\ $
with equality only if $X\cong \pn 1$.  Now in general
let $A\in |L|$.
Then by  using the exact sequence
$\ 0\to (s-1)L\to sL\to sL_A\to 0\ $
for $1\le s\le t$ we see that
$\displaystyle h^0(tL)\le \sum_{j=0}^t h^0(jL_A)$.
Thus by induction we have
$\displaystyle h^0(tL)\le \sum_{j=0}^t{j+n-2\choose
n-2}\frac{dj+n-1}{n-1}={t+n-1\choose
n-1}\frac{dt+n}{n}$.
\qed

Now assume that $L$ is very ample. Then we also have the following lower bound
\begin{equation}\label{easylower}
h^0(tL)\geq {{t+n+1}\choose{n+1}}\;\;{\rm for}\;\;t<d:=L^n.
\end{equation}
To see this note that we can assume $h^0(L)\geq n+2$ since otherwise
$(X,L)\cong(\pn
n,\sO_{\pn n}(1))$ and the assertion is clearly true. Then $X$
is embedded by
$|L|$ in $\pn {n+r}$ with $r>0$, so that we can project $X$ generically
one-to-one into $\pn
{n+1}$. Now, for any positive integer $t$,
$$h^0(tL)=h^0(\sO_{\pn {n+r}}(t)_{X})\geq h^0(\sO_{\pn {n+1}}(t)_{X'}),$$
where $X'$ is the image of $X$ in $\pn {n+1}$.
But if $t< d:=\deg_L(X)$ then  $h^0(\sO_{\pn {n+1}}(t)_{X})\geq
h^0(\sO_{\pn {n+1}}(t))$ since
the kernel of the restriction map has dimension $h^0(\sO_{\pn
{n+1}}(t-d))=0$. Thus we get
$h^0(tL)\geq h^0(\sO_{\pn {n+1}}(t))$, which is the bound as in
(\ref{easylower}).

Following Harris' presentation \cite{Harris} of Castelnuovo theory we can
significantly improve the above lower bound. Let us fix some notation. Let
$X_{n-i}$ be
the $(n-i)$-dimensional subvariety of $X$ obtained as transversal
intersection of $X$
with a general $\pn {n+r-i}$, $0\le i\le n$, and in particular
$X_n=X$.
Let, for $0\le i\le n$,
$$h_{X_{n-i}}(t) :=\dim({\rm Im}(H^0(\pn {n+r-i},\sO_{\pn {n+r-i}}(t))\to
H^0(X_{n-i},\sO_{X_{n-i}}(tL))),$$
$h^0(t):=h_{X_0}(t)$. By \cite[Lemma (3.1)]{Harris} one has, for a given
integer
$t\geq 0$,
$$h_X(t)\geq h_{X_{n-1}}(t)+h_X(t-1).$$
Iterating on $t$ we get  
$\ \displaystyle h_{X_{n-i}}(j)\geq \sum_{k=0}^{j}h_{X_{n-i-1}}(k)\ $
for $0\le i\le n$.
Iterating on $n$, we get
\begin{equation}\label{basicineq}
h^0(tL)\geq
h_X(t)\geq\sum_{k_{n-1}=0}^t\cdots\sum_{k_1=0}^{k_2}\sum_{k=0}^{k_1}h_{X_0}(k).
\end{equation}
Castelnuovo theory (see \cite[p. 94]{Harris}) gives
\begin{equation}\label{graphineq}
h_{X_0}(k)\geq h(k) := \min\{d,kr+1\}.
\end{equation}
(Note that the formula in \cite[p. 94]{Harris} is for a curve in $\pn {r}$,
whereas we are considering a curve in $\pn {r+1}$.)
Let $c:=\left[\frac{d-1}{r}\right]$, the integral part of $\frac{d-1}{r}$,
and let $R := d - 1 -cr$. Then the
graph of the function $h(k)$ looks like

\bigskip\bigskip\bigskip\bigskip

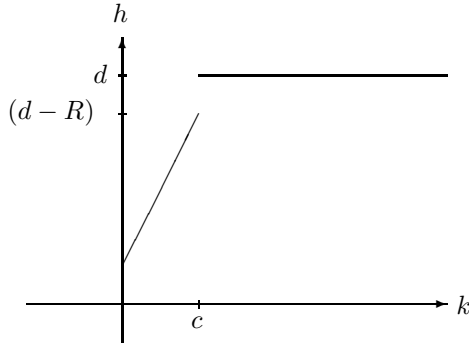
\begin{equation}\label{D1}
\setlength{\unitlength}{.01in}
\centerline{
\begin{picture}(230,100)(-50,-20)
\put (-50,0){\vector(1,0){220}}
\put(0,-20){\vector(0,1){160}}
\put(175,-5){$k$}
\put(-5,148){$h$}
\put(0,20){\line(1,2){40}}
\put(40,120){\line(1,0){130}}
\put(40,-2){\line(0,1){4}}
\put(36,-12){$c$}
\put(-2,120){\line(1,0){4}}
\put(-15,116){$d$}
\put(-2,100){\line(1,0){4}}
\put(-60,97){$(d-R)$}
\end{picture}
}
\end{equation}
where the oblique line is the graph of the equation $h=kr+1$.

\begin{lemma}\label{CastelL} Let $X$ be an $n$-dimensional irreducible
nondegenerate subvariety of $\pn {n+r}$. Let $L:=\sO_{\pn {n+r}}(1)_{X}$,
and let $d=L^n$ be the
degree of $X$
in $\pn {n+r}$. Let $c:=\left[\frac{d-1}{r}\right]$. Let $R$ be
the remainder defined as $R:=d-1-cr$. Then for any integer
$t\geq 0$ we have the lower bound
\begin{eqnarray}
h^0(tL)&\geq & r{{n+t}\choose{n+1}}+{{n+t}\choose
n}-r{{n+t-c-1}\choose{n+1}}+(R-r){{n+t-c-1}\choose n}\\ \nonumber
&=&\frac{tr+n+1}{t+n+1}{{t+n+1}\choose{n+1}}-r{{t+n-c}\choose{n+1}}+R{{t+n-c-1}
\choose n}. \nonumber
\end{eqnarray}
\end{lemma}
\proof Referring to the inequality in formula (\ref{graphineq}), we get for
any positive integer $a$
$$\sum_{k=0}^a h_{X_0}(k) \ge \sum_{k=0}^a h(k) = \sum_{k=0}^a kr+1 -
\sum_{k=c+1}^a (kr + 1 -d)$$
(see diagram (\ref{D1})).  Iterating the summation as in
formula (\ref{basicineq}) gives
\begin{eqnarray}
h^0(tL)&\geq &\sum_{k_{n-1}=0}^t\cdots\sum_{k_1 =
0}^{k_2}\sum_{k=0}^{k_1}(kr+1)
 - \sum_{k_{n-1}=c+1}^t\cdots\sum_{k_1=c+1}^{k_2}\sum_{k=c+1}^{k_1}(kr+1-d)
\nonumber \\
&=& \sum_{k_{n-1}=0}^t\cdots\sum_{k_1=0}^{k_2}\sum_{k=0}^{k_1}(kr+1) -
\sum_{j_{n-1}=0}^{t-c-1}\cdots\sum_{j_1=0}^{j_2}\sum_{j=0}^{j_1}((j+c+1)r+1-d)
\nonumber \\
&=& r \sum_{k_{n-1}=0}^t\cdots\sum_{k_1 = 0}^{k_2}\sum_{k=0}^{k_1} k
        +\sum_{k_{n-1}=0}^t\cdots\sum_{k_1 = 0}^{k_2}\sum_{k=0}^{k_1} 1
\nonumber \\
 & & \qquad - r
\sum_{j_{n-1}=0}^{t-c-1}\cdots\sum_{j_1=0}^{j_2}\sum_{j=0}^{j_1} j
-\sum_{j_{n-1}=0}^{t-c-1}\cdots\sum_{j_1=0}^{j_2}
\sum_{j=0}^{j_1}((c+1)r + 1 - d)\nonumber
\end{eqnarray}
By repeatedly using the combinatorial identity
$\ \displaystyle\sum_{i=0}^b {{i + m}\choose q} = {{b+m+1}\choose{q+1}} - 
{m\choose{q+1}}\ \ $
for any positive integers $b$, $m$, and $q$, with the usual convention that
${s\choose t} = 0$
whenever $t > s$, we get
\begin{eqnarray}\nonumber
h^0(tL)&\geq & r{{n+t}\choose{n+1}}+{{n+t}\choose
n}-r{{n+t-c-1}\choose{n+1}}+(d-1-r(c+1))
{{n+t-c-1}\choose n}\\ \nonumber
&=& r{{n+t}\choose{n+1}}+{{n+t}\choose n}-r{{n+t-c-1}\choose{n+1}}+(R-r)
{{n+t-c-1}\choose n} \\ \nonumber
&=&
\frac{tr+n+1}{t+n+1}{{t+n+1}\choose{n+1}}-r{{t+n-c}
\choose{n+1}}+R{{t+n-c-1}\choose n}.
\end{eqnarray}
\qed

\begin{rem*}\label{Rckod} It is easy to see that the right-hand side of
this last inequality is minimized when
$R=0$ and $c=1$, and therefore the bound in Lemma (\ref{CastelL}) yields the
simpler form
\begin{equation}\label{simpler}
h^0(tL)\geq \frac{tr+n+1}{t+n+1}{{t+n+1}\choose{n+1}}-r{{t+n-1}\choose{n+1}}.
\end{equation}
Note also that if $X$ has nonnegative Kodaira dimension, then $d\geq rn+2$
(see e.g.,
\cite[(8.1.3)]{Book}). Thus $\frac{d-1}{r}>n$, so that
$c=\left[\frac{d-1}{r}\right]\geq
n$. Therefore in this case we can use the bound in (\ref{CastelL}), with
$R=0$, $c=n$, in
the form
$$h^0(tL)\geq \frac{tr+n+1}{t+n+1}{{t+n+1}\choose{n+1}}-r{t\choose {n+1}}.$$
\end{rem*}

We have the following general fact.

\begin{prop}\label{CastelGen} Let $X$ be a nondegenerate  irreducible
$n$-dimensional
subvariety of $\pn {n+r}$. Let $L:=\sO_{\pn {n+r}}(1)_{X}$. Let
$Y$ be a $k$-dimensional irreducible subvariety of $X$ of degree
$\delta:=L^k\cdot Y$. Assume that $k>0$. Let $t> 0$ be an integer such that
$$\frac{tr+n+1}{t+n+1}{{t+n+1}\choose{n+1}}-r{{t+n-1}\choose{
n+1}}>
\frac{\delta t+k}{t+k}{{t+k}\choose k}.$$
Then $tL\otimes \sJ_Y$ has a section not identically zero  on $X$.\end{prop}
\proof Let $L_Y$ be the restriction of $L$ to $Y$. By Lemma (\ref{upper})
applied to $L_Y$ we have
\begin{equation}\label{ub}
h^0(tL_Y)\leq\frac{\delta t+k}{t+k}{{t+k}\choose k}.
\end{equation}
Suppose that $h^0(tL\otimes\sJ_Y)=0$. Then $h^0(tL)\leq h^0(tL_Y)$. Thus by
combining the inequalities (\ref{ub}) and (\ref{simpler}) we get
$$\frac{\delta t+k}{t+k}{{t+k}\choose k}\geq
\frac{tr+n+1}{t+n+1}{{t+n+1}\choose{n+1}}-r{{t+n-1}\choose{n+1}},$$
contrary to the inequality assumed in the proposition.
\qed

Let us now make explicit the bound in (\ref{CastelGen}) in the case when
$Y$ is a divisor on $X$.

\begin{prop}\label{CastelP} Let $X$ be a nondegenerate  irreducible
$n$-dimensional
subvariety of $\pn {n+r}$. Let $L:=\sO_{\pn {n+r}}(1)_{X}$. Let $D$ be an
irreducible divisor
of degree
$\delta:=L^{n-1}\cdot D>1$. Thus, for $t\geq 1$, the inequality
$\ \displaystyle t>\frac{n}{r+1}(\delta-1)-n+1\ $
implies that $tL-D$ has a section not identically zero on $X$.\end{prop}
\proof The inequality in (\ref{CastelGen}) becomes, in case $k=n-1$,
\begin{equation}\label{uplozero}
\frac{tr+n+1}{t+n+1}{{t+n+1}\choose{n+1}}-r{{t+n-1}\choose{n+1}}>\frac{\delta
t+n-1}{t+n-1}{{t+n-1}\choose{n-1}}.
\end{equation}
Now by a simple calculation
(\ref{uplozero}) gives
$-t-rt-2n+\delta n+1-rn+r< 0$,
or, solving in $t$ and simplifying,
$\displaystyle t>\frac{n}{r+1}(\delta-1)-n+1.$ \qed

The following example shows that Proposition (\ref{CastelP}) is
sharp.

\begin{example}\label{ExP}  Let $X:=\pn2\times \pn2$, with $L := \sO(1,1)$,
and choose $D \in |\sO(2,0)|$.
Then  $tL - D \approx \sO(t-2, t)$ has a non-trivial section if and only if
$t \geq 2$.
Embed $X$ in $\pn 8$ by the Segre mapping.  Then $n = r = 4$, and an easy
calculation gives
$\delta := L^3\cdot D = 6$.

In this case we see that the hypothesis of (\ref{CastelP}) is satisfied if
$t=2$, and $2L - D$ has a non-zero section; whereas
the hypothesis fails if $t=1$ and $L-D$ has no non-trivial sections. Thus
the inequality  in (\ref{CastelP}) cannot be weakened.
\end{example}

\begin{rem*}\label{cod2}We follow the notation and assumptions of
(\ref{CastelGen}). In general, it is   hard to make
the bound in (\ref{CastelGen}) explicit in $t$. If $Y$ has codimension two
in $X$,
then, after simplification,  the
condition for $h^0\left((\delta-1)L\otimes\sJ_Y\right)$ to be positive  becomes
$$\ -rn^2-3n^2+\delta n^2-2rn\delta+9n-4\delta
n+5rn-6-6r+5\delta+5r\delta-r\delta^2-\delta^2 > 0.$$
\end{rem*}

Let $X\subset \pn N$ be a nondegenerate smooth connected $n$-fold.
Let $\deg(X)=d$ and denote by $L$ the restriction of $\sO_{\pn N}(1)$
to $X$. From Lemma (\ref{Mumford}) we know that $\sO_{\pn
N}(d)\otimes\sJ_X$ is spanned by global sections.

\noindent {\bf Problem.} What can we say about the smallest
integer $t>0$ such that $h^0(\sO_{\pn N}(t)\otimes\sJ_X)>0$?

We define the {\em lower degree} in $\pn N$, $\delta_N$, of a
subvariety $X\subset \pn N$ to be the smallest positive integer $t$ such
that $h^0(\sO_{\pn N}(t)\otimes\sJ_X)>0$. One consequence of the
above results  is that under modest conditions there must be some form
of much lower degree than $d$ vanishing on $X$.

\begin{prop}\label{lowerbound}Let $X\subset \pn N$ be a nondegenerate smooth
 connected $n$-fold of degree $d$. Let $L=\sO_{\pn N}(1)_{X}$ and let
$\delta_N$ be the degree of the lowest-degree homogeneous form vanishing on
$X$. Then
$\delta_N$ satisfies the inequality
$\ \displaystyle
\frac{1}{\delta_N-1}\left(\frac{(\delta_N+N-1)\cdots(\delta_N+n-1)}
{N\cdots(n+1)}-n\right)\leq d.$
\end{prop}
\proof If $t$ is a positive integer for which $h^0(\pnsheaf
Nt\otimes\sJ_X)=0$ then
$${{N+t}\choose N} = h^0(\sO_{\pn N}(t)) \le h^0(X, tL).$$
By applying Lemma (\ref{upper}) we get
$\ \displaystyle {{N+t}\choose N} \le {{td+n}\over{t+n}}{{t+n}\choose n}$.
An easy calculation shows that this inequality is equivalent to
$\ \displaystyle{1\over t}\left({{(t+N)\cdots(t+n)}\over{N\cdots(n+1)}} -
n\right) \le d$.
 From the definition of $\delta_N$ it follows that
$h^0(\pnsheaf N{\delta_N-1}\otimes\sJ_X)=0$,
and substituting $t = \delta_N - 1$ in the last inequality completes the
proof. \qed

For surfaces here is the explicit bound.

\begin{corollary}\label{lbsurf}Let $X$ be a nondegenerate smooth
 connected surface of  degree $d$ in $\pn N$. Assume $N\geq 5$.
Let $\delta_N$ be the degree of the lowest-degree homogeneous form
vanishing on $X$.
 Then $\delta_N^3+11\delta_N^2+46\delta_N+96\leq 60d$.
\end{corollary}
\proof We apply the bound in (\ref{lowerbound}) with $n=2$. Since $N \ge
5$, we get
$$\frac{1}{\delta_N-1}\left(\frac{(\delta_N+4)\cdots(\delta_N+1)}
{60}-2\right)\leq d.$$
After simplifying this becomes
$\ \displaystyle (\delta_N +4)(\delta_N + 3)(\delta_N + 2)(\delta_N + 1) - 120 \le
60d(\delta_N - 1)$.
Since $\delta_N\geq 2$ we can divide both sides by $\delta_N - 1$ to obtain
the desired
result.\qed

For example, if $X$ is  of degree $21$ then $\delta_N\leq 6$ and hence
there is a form of the sixth degree vanishing on $X$. Or again, if $X$ is of 
degree $10$,$000$ there is a form of degree $80$ vanishing on
$X$.

Some other special cases are as follows. For threefolds in $\pn N$ with
$N\geq 5$ the corresponding bound as in Corollary (\ref{lbsurf}) is
$\ \delta_N^2+10\delta_N+36\leq 20d.$
For threefolds in $\pn N$ with $N\geq 6$ the bound becomes
$\ \delta_N^3+15\delta_N^2+86\delta_N+240\leq 120d.$

\section{Some general structure results for
projections}\label{StructureSec}
\addtocounter{subsection}{1}\setcounter{theorem}{0}

In this section we discuss some general properties of projections from
a $k$-dimensional subvariety $Y$ of a given polarized variety $X$. We always assume that
$k>0$.
In
\S \ref{Divisorial} and \S \ref{Linear} we will present some more refined
results in the cases when $Y$ is either  a divisor or a linear $\pn k$.

\begin{prgrph*}{General set-up of morphisms}\label{SetUp}
Let $L$ be a very ample line bundle on $X$, a smooth connected variety
of dimension $n\geq 2$.  Let $Y$
be a $k$-dimensional connected submanifold of $X$. {\em We always assume that 
$k>0$.} We will denote by
$\sJ_Y$ the ideal sheaf of $Y$ in $X$.

Let $\grs:\oline{X}\to X$ be the blowing up of $X$ along $Y$ and set
$E=\grs^{-1}(Y)$. Let $\psi:\oline{X}\to \psi(\oline X)$ be the surjective
rational map given
by $|\grs^*L-E|$. We refer to the mapping $\psi$ as the {\em
projection from $Y$ associated to $L$}. If $L\otimes\sJ_Y$ is spanned
by its global sections, then
$\grs^*L-E$ is spanned on
$\oline{X}$ and $\psi$ is a morphism. We have the Remmert-Stein factorization
$\psi= s\circ \phi$ of $\psi : \oline X\to \psi(\oline X)$, where
$\phi:\oline X\to Z$ is a morphism with
connected  fibers onto a normal variety $Z$ and $s:Z\to \psi(\oline X)$ is
a finite morphism. We will refer to $\phi:\overline{X}\to Z$ as the morphism 
associated to $L\otimes\sJ_Y$.
Note there is an ample and spanned line bundle $\sH$ on $Z$
such that $\grs^*L-E\approx\phi^*\sH$.
We have the following commutative diagram
$$
\begin{array}{ccccc}
E&\hookrightarrow&\!\!\!\!\!\!\!\!\!\!\oline{X}&&\\
\downarrow&&\!\!\!\!\!\!\downarrow\grs&
\!\!\!\!\!\!\!\!\!\!\!\!\!\!\!\stackrel{\phi}{\searrow}&\\
Y&\hookrightarrow&X\stackrel{\vphi}{\longrightarrow} &Z
\end{array}
$$
where $\vphi$ is the connected part of the rational mapping associated to
$|L\otimes\sJ_Y|$.
\end{prgrph*}

We need the following technical lemmas.

\begin{lemma}\label{differentials} Let  $Y$
be a connected $k$-dimensional submanifold of $X$, a connected projective
manifold of
dimension $n\geq 2$. Assume that $k>0$. Let $N:=N_{Y/X}$ be
the normal bundle of $Y$ in $X$. Let $L$ be a line bundle on $X$.
Assume that $L\otimes\sJ_Y$ is spanned by a vector subspace, $V$, of
$\Gamma(L\otimes\sJ_Y)$  and $c_{n-k}(N^*(L))=0$. Then a general element
$D\in |V|$ is smooth.\end{lemma}
\proof A general $D\in |V|$ is smooth on $X\setminus Y$ by
Bertini's theorem.

A given $D\in |L\otimes\sJ_Y|$ is smooth at a point $y\in Y$ if the
differential in local coordinates of the defining equation of $D$ is not
zero at  $y\in Y$. From the exact sequence
$$0\to \sJ_Y^2\otimes L\to\sJ_Y\otimes
L\stackrel{\partial}{\to}N^*(L)\to 0$$ we see that $D$ is smooth on
$Y$ if the image $\partial(s)$ of $s$ defining $D$ in $N^*(L)$ is
nowhere zero. Since $c_{n-k}(N^*(L))=0$, a general $s\in
V$ goes to a nowhere vanishing section
$\partial(s)$ in $N^*(L)$.\qed

\begin{lemma}\label{imdiv}Let  $Y$
be a connected $k$-dimensional submanifold of $X$, a connected projective
manifold of
dimension $n\geq 2$. Assume that $k>0$. Let $L$ be a very
ample line bundle on $X$. Assume that $L\otimes\sJ_Y$ is spanned by
global sections and
let $D$ be a smooth element of $|L\otimes\sJ_Y|$. Let
$\grs:\oline{X}\to X$ be the blowing up of $X$ along $Y$ and let
$\oline{D}$ be the proper transform of $D$ under $\grs$. Let
$\phi:\oline{X}\to Z$ be the morphism associated to $L\otimes\sJ_Y$ as
in {\rm (\ref{SetUp})}. Let $\phi_{\oline{D}}$ be the restriction of $\phi$
to $\oline{D}$. Then $\phi_{\oline{D}}$ has lower dimensional image if
and only if $\phi$ has lower dimensional image.
\end{lemma}
\proof  We follow the  notation from (\ref{SetUp}).
Assume that $\phi$ has lower dimensional image. We have
$\oline{D}\in|\grs^*L-E|$, $E=\grs^{-1}(Y)$. Note that $\oline{D}$ is
the pullback of some divisor $\sD\subset Z$. From
$\ \displaystyle \dim\phi(\oline{D})<\dim\phi(\oline{X})<\dim X\ $
we get $\dim\phi(\oline{D})<\dim\oline{D}$.

To show the converse, note that by definition of $\oline{D}$ one has
$\dim\phi(\oline{D})=\dim\phi(\oline{X})-1$. Thus the assumption
$\dim\phi(\oline{D})<\dim\oline{D}=n-1$ gives the result.\qed

Let us note some further general properties of the morphism
$\phi$. The notation is as in (\ref{SetUp}).
\begin{enumerate}
\item[1.] (Divisorial case) If $Y$ is a divisor, then $X\cong\oline{X}$.
\item[2.] (Linear case) Assume that $(Y,L_Y)\cong\pnpair k 1$ and that 
$\Gamma(L)$ embeds $X$ in $\pn {n+r}$. Since $Y$ is a linear space it
follows that $\displaystyle
L\otimes\sJ_Y\;\;{\rm is\; spanned\; by\; global\; sections}$.
The mapping given by $\Gamma(L\otimes\sJ_Y)$ coincides off of $Y$ with
the restriction to $X$
  of  the projection of $\pn {n+r}$ to $\pn {n+r-k-1}$ from $Y$.
\item[3.] (Smooth case) If $Y$ is smooth then $\delta
L\otimes\sJ_Y$, $\delta=L^k\cdot Y$, is spanned by global sections by
Lemma (\ref{Mumford}).
\end{enumerate}

We have the following crude structure theorem in the case when the
projection has lower dimensional image.

\begin{theorem}\label{Structure} Let  $Y$
be a connected $k$-dimensional submanifold of $X$,
a connected projective manifold of dimension $n\geq 2$. 
 Assume that $k>0$. Let $L$ be a very
ample line bundle on $X$. Assume that $L\otimes\sJ_Y$ is spanned by
its global sections.  Let $\grs:\oline{X}\to X$ be the blowing up of $X$
along $Y$.  Let $\phi:\oline{X}\to Z$ be the morphism
associated to $L\otimes\sJ_Y$ as in {\rm (\ref{SetUp})}. Let
$E:=\grs^{-1}(Y)$ be the
exceptional divisor. Assume $n >\dim Z$. Then we have:
\begin{enumerate}
\em\item\em $\phi(E) = Z$;
\em\item\em $Z$ is uniruled if $\cod_XY>1$;
\em\item\em $Z$ is unirational if $Y$ is unirational;
\em\item\em the restriction, $E_F$, of $E$ to any fiber $F$ of $\phi$
is an ample divisor on $F$.
\end{enumerate}\end{theorem}
\proof To show $1$), assume by contradiction that the restriction,
$\phi_E:E\to Z$, of $\phi$ to $E$ is not surjective. Take a point $x\in
Z\setminus\phi(E)$ and let $F_x=\phi^{-1}(x)$ be the fiber on $x$. Then
the restriction $(\grs^*L-E)_{F_x}$ is trivial. But $E_{F_x}\cong
\sO_x$ since $x\notin \phi(E)$, so that $(\grs^*L-E)_{F_x}\cong
(\grs^*L)_{F_x}\cong L_{\grs(F_x)}$ is ample, where the last isomorphism
follows from the fact that $F_x$ goes isomorphically to $X$ under
$\grs$, since $F_x\cap E=\emptyset$. Thus $(\grs^*L-E)_{F_x}$ is both
ample and trivial, an absurdity that
contradicts the assumption $n>\dim Z$.

To show $2$), note that since $\cod_XY>1$, $\grs$ is not an isomorphism
and the exceptional divisor is uniruled. This means that there exists
an $(n-2)$-dimensional variety $V$ and a rational map $V\times \pn 1\to
E$ which is dominant. Since $\phi(E)=Z$ by $1$), we get a
dominant map $V\times \pn 1\to Z$, i.e., $Z$ is uniruled.

To show $3$), recall that $E$ is  birational to $Y\times \pn {n-k-1}$.
Since $Y$ is unirational we have a dominant rational map
$\pn k\to Y$. Therefore, combining with the surjective map $\phi_E:E\to Z$,
we get a dominant rational map $\pn k\times\pn
{n-k-1}\to Z$. This implies that $Z$ is unirational.

To show $4$), take a fiber $F$ of $\phi:\oline{X}\to Z$. If the
restriction $\grs_F$ of the blowing up map is finite-to-one then
$\grs^*_FL$ is ample and the assertion is clear. It is easy to see that
$\grs_F$ is finite. If not, then it follows that there is a positive
dimensional fiber, $f$, of $\grs_F:F\to \grs(F)$. This implies that $f$
is contained in a fiber of $E\to Y$. But $\grs^*L-E$ is ample on
fibers of $E\to Y$. On the other hand, since $f\subset F$, the line
bundle $\grs^*L-E$ is trivial on $f$.\qed

We need the following lemma.

\begin{lemma}\label{Equivalence} Let  $Y$ be a connected $k$-dimensional
submanifold of $X$, a connected projective manifold of
dimension $n\geq 2$. Assume that $k>0$. Let $L$ be a very ample
line bundle on $X$. Let $\sJ_Y$ be the ideal sheaf of $Y$ in $X$ and let
$N:=N_{Y/X}$ be the normal bundle of $Y$ in $X$. Then  $L\otimes\sJ_Y$
is spanned and $N^*(L)$ is trivial if and only if $Y$ is the complete
intersection of $n-k$ divisors $D_1,\ldots,D_{n-k}\in |L|$.
\end{lemma}
\proof The ``if'' part is straightforward. As to the converse, consider
the exact sequence
$$0\to L\otimes\sJ_Y^2
\to L\otimes\sJ_Y\stackrel{\partial}{\to}N^*(L)\to 0.$$
Set $w:=n-k$. Since $N^*(L)=\oplus^w\sO_Y$ and $L\otimes\sJ_Y$ is
spanned we can find sections $s_1,\ldots,s_w\in \Gamma(L\otimes\sJ_Y)$
defining $w$ divisors $D_1,\ldots,D_w\in |L\otimes\sJ_Y|$ on $X$
containing $Y$.

For each $i=1,\ldots,w$, $D_i$ is smooth on $X\setminus Y$ by Bertini's
theorem.

For each $i=1,\ldots,w$, $D_i$ is smooth at a point $y\in Y$ if the
differential in  local coordinates of the defining equation of $D_i$ is
not zero at $y$. From the exact sequence above we see that $D_i$ is
smooth on $Y$ if the image $\partial(s_i)$ of $s_i$ defining $D_i$ in
$N^*(L)$ is nowhere zero, $1\le i\le w$. Since $N^*(L)$ is trivial
we can find the sections $s_1,\ldots,s_w\in \Gamma(L\otimes\sJ_Y)$
such that $\partial(s_1),\ldots,\partial(s_w)$ are independent in
$N^*(L)$. It follows that $D_1,\ldots,D_w$ are smooth as well as the
intersection $D_1\cap\ldots\cap D_w$ is smooth.

Since $\dim X\geq 2$, from the exact sequence
$\displaystyle 0\to -L\to\sO_X\to\sO_{D_i}\to 0$,
we see that $h^0(\sO_{D_i})=1$, so  the $D_i$'s are connected,
$1\le i\le w$.

Since $D_1\cap\ldots\cap D_w$ is at least one-dimensional, we know by 
the Lefschetz hyperplane section theorem that $D_1\cap\ldots\cap D_w$ 
is connected. Since it is also smooth of dimension $\dim Y$ and contains
 $Y$ we  conclude that $Y$ is the complete intersection of
$D_1,\ldots,D_w$.\qed

We can prove now the following more refined structure result, which
gives a general lower bound for the dimension of the image of the
projection.

\begin{theorem}\label{Dimension} Let  $Y$ be a connected $k$-dimensional
submanifold of $X$, a connected projective manifold of
dimension $n\geq 2$.  Assume that $k>0$. Let $L$ be a very
ample line bundle on $X$. Assume that $L\otimes\sJ_Y$ is spanned by
its global sections.  Let $\grs:\oline{X}\to X$ be the blowing up of $X$
along $Y$.  Let $\phi:\oline{X}\to Z$ be the morphism  from
associated to $L\otimes\sJ_Y$ as in {\rm (\ref{SetUp})}. Then $\dim Z\ge
n-k-1$, with equality   if and only
if  $Z\cong\pn {n-k-1}$ and $Y$ is the complete intersection of $n-k$
divisors $D_1,\ldots,D_{n-k}\in |L|$.
\end{theorem}
\proof Set $w:=n-k$. Let $E:=\grs^{-1}(Y)$
be the exceptional divisor.
Set  $\oline{L}:=\grs^*L$ and $N:=N_{Y/X}$. Recall that $E\cong{\Bbb
P}(N^*(L))$. Thus $\pi:=\grs_E:E\to Y$ is a $\pn {w-1}$-bundle. Let
$\xi$ be the tautological line bundle of ${\Bbb P}(N^*(L))$. Notice
that $\xi\cong(\oline{L}-E)_E$. Since $(\oline{L}-E)_{\pn
{w-1}}\cong\sO_{\pn
{w-1}}(1)$ it follows that each fiber $\pn {w-1}$ of $\pi:E\to Y$ maps
isomorphically  under the map $\psi:\oline{X}\to \psi(\overline{X})$ given
by $|\oline{L}-E|$, and hence maps isomorphically into $Z$ under the
morphism $\vphi$ associated to $L\otimes\sJ_Y$. 
 This shows that $\dim Z\geq w-1$. If $\dim Z=w-1$, it
follows that $Z \cong \pn{w-1}$.

It also follows that $Y$ is a complete intersection. For we have a
surjective map of locally free sheaves
$\ \displaystyle\oplus^{\dim Z+1}\sO_{\oline{X}}\to \oline{L}-E\to 0$,
and, restricting to $E$, we have a surjection
$\ \displaystyle\oplus^{\dim Z+1}\sO_E\to (\oline{L}-E)_E\to 0$.
Consider the $\pn {w-1}$-bundle $\pi:E\to Y$. Notice that
$(\oline{L}-E)_E\cong\xi$, the tautological line bundle of
$E\cong{\Bbb P}(N^*(L))$. By pushing forward under $\pi$, we get a
surjection
$$\grb:\oplus^{\dim Z+1}\sO_Y\to N^*(L)=\pi_*\xi\to 0.$$ By comparing
the ranks, since $N^*(L)$ has rank ${\rm cod}_XY=w=\dim Z+1$, we
conclude that $\grb$ is an isomorphism, i.e., $N^*(L)$ is the trivial
bundle. Thus, since $L\otimes\sJ_Y$ is spanned by global sections,
Lemma (\ref{Equivalence}) applies to give the result.

Next, we show that if $Y:=D_1\cap\ldots\cap D_w$ is the
complete intersection of $w$ divisors $D_1,\ldots,D_w\in |L|$,
then $\dim Z = w-1$. We first observe that
$N\cong\oplus^w{L_Y}$, so that $N^*(L)\cong\oplus^w\sO_Y$ is
trivial. Consider the exact sequence
$\ 0\to L\otimes\sJ_Y^2\to L\otimes\sJ_Y\to N^*(L)\to 0.\ $
Since the morphism $L\otimes\sJ_Y\to N^*(L)$ is surjective at the
sheaf level, $L\otimes\sJ_Y$ is spanned by global sections and
$N^*(L)$ is trivial, it follows that the induced map
$\ \displaystyle\gra:\Gamma(L\otimes\sJ_Y)\to\Gamma(N^*(L))\to 0\ $
is surjective.

For any integer $m\geq 2$, consider the exact sequence
$$0\to L\otimes\sJ_Y^{m+1}\to L\otimes\sJ_Y^m\to S^m(N^*)\otimes L\to
0.$$
Since $S^m(N^*)\otimes L\cong\oplus{L_Y}^{-(m-1)}$ we have
$h^0(S^m(N^*)\otimes L)=0$, for $m\geq 2$,
and therefore we get
$$\Gamma(L\otimes\sJ_Y^2)\cong\cdots\cong
\Gamma(L\otimes\sJ_Y^m),\;\;m\geq 2.$$
If $h^0(L\otimes\sJ_Y^2)\neq 0$ we thus find a section of $L$
vanishing on $Y$ of any given order $m\geq 2$, which is absurd.
Therefore we conclude that $h^0(L\otimes\sJ_Y^2)=0$ and hence
$\Gamma(L\otimes\sJ_Y)$ injects in $\Gamma(N^*(L))$, i.e., the map
$\gra$ is an isomorphism. Thus
$h^0(L\otimes\sJ_Y)=h^0(\oline{L}-E)=w$. Since $\oline{L}-E$ is
spanned and gives the projection $\psi:\oline{X}\to \psi(\overline{X})$, we
thus conclude that the $\dim Z = w-1$.
\qed

There are many results from adjunction theory \cite{Book}
describing all varieties with a given hyperplane section. Combining
these results with Theorem (\ref{Dimension}) gives many consequences.
By way of illustration we give two useful corollaries.

\begin{corollary}\label{Pk} Let  $Y$ be a connected $k$-dimensional
submanifold of $X$, a connected projective manifold of
dimension $n\geq 2$.  Let $L$ be a very
ample line bundle on $X$. Assume that $k>0$ and $Y$ is a linear $\pn k$
with respect to $L$, i.e., $L^k\cdot Y=1$. Let $\grs:\oline{X}\to X$ be the
blowing up of $X$ along $Y$. Let
$\phi:\oline{X}\to Z$ be the morphism  associated to $L\otimes\sJ_Y$ as in
{\rm (\ref{SetUp})}. Then $\dim Z=n-k-1$ if and only if $(X,L)\cong(\pn
n,\sO_{\pn n}(1))$.\end{corollary}
\proof Assume $\dim Z=n-k-1$. Then, by (\ref{Dimension}), $2$), $Y$ is
the complete intersection of $n-k$ divisors $D_1,\ldots,D_{n-k}\in
|L|$. Since $Y$ is a linear $\pn k$ it thus follows that $(X,L)\cong
(\pn n,\sO_{\pn n}(1))$.

If $(X,L)\cong(\pn n,\sO_{\pn n}(1))$, the projection from $Y=\pn k$
has an $(n-k-1)$-dimensional image.\qed

\begin{corollary}\label{KP1} Let  $Y$ be a connected $k$-dimensional
submanifold of $X$, a connected projective manifold of
dimension $n\geq 3$. Assume that $k>0$. Let $L$ be a very ample line bundle 
on $X$. Assume that $L\otimes\sJ_Y$ is spanned by its global sections.
Let $\grs:\oline{X}\to X$ be the blowing up of $X$ along $Y$.
Let $\phi:\oline{X}\to Z$ be the morphism associated to $L\otimes\sJ_Y$
as in {\rm (\ref{SetUp})}.  Assume that $Y$ is a $K(\pi,1)$ and $k\geq 2$.
Then $\dim Z\geq n-k$.\end{corollary}
\proof By (\ref{Dimension}) either we are done or $Y$ is a complete
intersection of $n-k$ divisors $D_1,\ldots,D_{n-k}\in |L|$. Since $Y$
is a $K(\pi,1)$ with $\dim Y\geq 2$, this is not possible by a result of
the fourth author \cite{So2}.\qed

Under special conditions on the cohomology of $Y$, we get stronger lower
bounds for the image dimension of the projection. We restrict our attention
to the case in which
$\dim Z \geq n-k$, since the case   $\dim Z = n-k-1$ was covered in
(\ref{Dimension}).

\begin{prop}\label{BetterB} Let  $Y$ be a connected $k$-dimensional
submanifold of $X$, a connected projective manifold of
dimension $n\geq 2$. Assume that $k>0$. Let $L$ be a very
ample line bundle on $X$. Assume that $L\otimes\sJ_Y$ is spanned by
its global sections.  Let $\grs:\oline{X}\to X$ be the blowing up of $X$
along $Y$.  Let $\phi:\oline{X}\to Z$ be the morphism
associated to $L\otimes\sJ_Y$ as in {\rm (\ref{SetUp})}. Assume that
$h^2(Y,\rat)_{\rm alg}=1$ {\rm (}or equivalently
 ${\rm Pic}(Y)\otimes\rat\cong\rat${\rm )}. If $\dim Z \geq n-k$, then
$\displaystyle\dim Z\geq n-k+\frac{k}{n-k}-1$.
In particular, if  $\dim Z={\rm cod}_XY=n-k$, then $n \ge 2k$.
\end{prop}
\proof Let $w:=n-k$ and $N:=N_{Y/X}$. As in the proof of
(\ref{Dimension}), we have a surjective vector bundle map
$\ \displaystyle \oplus^{\dim Z+1}\sO_Y\to N^*(L)\to 0.\ $
This gives a natural map
$\rho:Y\to{\rm Grass}(w,\dim Z+1)$ of $Y$ in the Grassmannian of the 
$w$-dimensional quotients of $\comp^{\dim Z+1}$. We claim that the map
$\rho$ is finite. Indeed, to see this, notice that ${\rm
det}(N^*(L))=\rho^*\sP$, where $\sP$ is an ample line bundle, the
Pl\"{u}cker bundle, on ${\rm Grass}(w,\dim Z+1)$.
Since $\rho$ is a not trivial map, ${\rm det}(N^*(L))$ is spanned
and not trivial. Since ${\rm Pic}(Y)\otimes\rat\cong\rat$, we thus
conclude that $\rho^*\sP$ is ample.
Let $F$ be a connected component of a positive dimensional fiber of
$\rho$. Then $(\rho^*\sP)_F\cong\sO_F$. This contradicts the ampleness of
$\rho^*\sP$. Thus
$$k=\dim Y\leq\dim{\rm Grass}(w,\dim Z+1)=w(\dim Z+1-w)$$
 gives the desired inequality.

If $\dim Z=n-k$, we have $k\leq n-k={\rm cod}_XY$, which is the same as
$2k\leq n$.
\qed

\begin{theorem}\label{ThmB}  Let  $Y$ be a connected $k$-dimensional
proper submanifold of $X$, a connected projective manifold of
dimension $n\geq 2$.  Assume that $k>0$.
Let $L$ be a very
ample line bundle on $X$. Assume that $L\otimes\sJ_Y$ is spanned by
its global sections.  Let $\grs:\oline{X}\to X$ be the blowing up of $X$
along $Y$.  Let $\phi:\oline{X}\to Z$ be the morphism
associated to $L\otimes\sJ_Y$ as in {\rm (\ref{SetUp})}. If
$h^{2j}(Y,\rat)_{\rm alg}=1$ for $j\leq w:=n-k$, then either $Z\cong\pn
{n-k-1}$  with $Y$ the
complete intersection of $n-k$ divisors in $|L|$  or $\dim Z\geq k$.
\end{theorem}
\proof By Theorem (\ref{Dimension}) we can assume that if the theorem
is false then
\begin{equation}\label{standingAssumptions}
k-1\ge \dim Z\geq w.
\end{equation}
 Let $E=\grs^{-1}(Y)$ be the exceptional divisor of
$\grs:\oline{X}\to X$ and let $\phi_E:E\to Z$ be the restriction to $E$
of the morphism $\phi:\oline{X}\to Z$. From (\ref{Structure}), $1$),
we know that $\phi_E$ is surjective. Let $N:=N_{Y/X}$ be the normal
bundle of $Y$ in $X$. Set $\oline{L}=\grs^*L$. Let $\xi$ be the
tautological line bundle of ${\Bbb P}(N^*(L))\cong E$. Notice that
$\xi\cong(\oline{L}-E)_E$.
Therefore for each fiber $\pn {w-1}$ of the $\pn {w-1}$-bundle $\pi:E\to
Y$ we have
$(\oline{L}-E)_{\pn {w-1}}\cong\sO_{\pn {w-1}}(1).$
This implies that each fiber $F$ of $\phi_E$ meets $\pn {w-1}$ in at
most one point. It thus follows that $F$ goes isomorphically to $\pi(F)$ under
$\pi$. Since $\xi_F\cong\sO_F$, we
get a surjective map
$(\pi^*N^*(L))_F\to\xi_F\cong\sO_F\to 0.$
Letting $F':=\pi(F)$, we have by the above $F\cong F'$ and therefore pushing
forward under $\pi$ we get a surjective map
\begin{equation}\label{Onto}
N^*(L)_{F'}\to \sO_{F'}\to 0.
\end{equation}
We claim that
\begin{equation}\label{Chern0}
c_w(N^*(L))=0.\end{equation}
To see this, let $F$ be a general fiber of $\phi_E$ and $F'=\pi(F)$.
 Note that $\dim F'=\dim E-\dim Z\ge n-1-(k-1)=w$. In view of this and
the assumption that
$h^{2j}(Y,\rat)_{\rm alg}=1$, $j\leq w$, we see that
it is enough to note that from (\ref{Onto}) it immediately follows that
$
c_w(N^*(L)_{F'})=0$.

Thus Lemma (\ref{differentials}) implies that there exists a smooth
divisor $D\in |L\otimes\sJ_Y|$. Since $\dim X\geq 2$, from the exact
sequence
$\ 0\to -L\to\sO_X\to\sO_D\to 0$, we see that $h^0(\sO_D)=1$, so $D$ is
connected. Let $H\in|\oline{L}-E|$ be the divisor corresponding to $D$.
I.e., $H=\phi^*Z_1$, where $Z_1\in |\sH|$ and $\oline{L}-E\approx \phi^*\sH$
for some ample and spanned line bundle $\sH$ on $Z$. Note that by the
generalized Seidenberg theorem (see e.g., \cite[(1.7.1)]{Book}), $Z_1$ is
irreducible
and normal since $Z$ is irreducible and normal. Notice also
that $\grs^*D\approx H+E$.

Let $X_1:=D$. By construction, $Y\subset X_1$. Furthermore
the blowing up $\grs:\oline{X}\to X$ induces a blowing up map
$\grs_1:\oline{X}_1\to X_1$ of $X_1$ along $Y$. We can also consider the
morphism, $\phi_1:\oline{X}_1\to Z_1$,  associated to
$L_{X_1}\otimes\sJ_Y$, where $L_{X_1}$ is the restriction $L_{X_1}$ of
$L$ to $X_1$. Note that $\phi_1$ is onto, $\dim X_1=n-1$, $\dim Z_1=\dim
Z-1$. Hence in particular $\dim Z_1<k$, i.e., (\ref{standingAssumptions})
is  preserved passing from $Z$ to $Z_1$.

Thus, starting from $X_1=D$, $Z_1$,
$\phi_1:\oline{X}_1\to Z_1$, $Y\subset X_1$, we proceed in such a way that
from the
initial data $$(k,\dim Z, w)$$
we reach, after $w-1$ steps, the data (recall that we are working under the
initial assumption that $\dim Z\geq w$)
$$(k,\dim Z- w+1,1).$$
I.e.,  $Y$ is a divisor in $X_{w-1}$ with
$X_{w-1}\cong\oline{X}_{w-1}$, and the morphism
$\phi_{w-1}:\oline{X}_{w-1}\to Z_{w-1}$ has image of dimension $\dim
Z-w+1$. In particular, since $X_{w-1}\cong\oline{X}_{w-1}$, we can
restrict $\phi_{w-1}$ to $Y$, so that we get a surjective map from $Y$ to
$Z_{w-1}$ (see (\ref{Structure}), $1$). By assumption we have that 
$h^2(Y,\rat)_{\rm alg}=1$, and therefore that $\dim Y=\dim Z_{w-1}$. Thus using
(\ref{standingAssumptions})
we have
$$ k=\dim Z_{w-1}=\dim Z-w+1=\dim Z - n+k+1$$
which gives that $n=\dim Z+1$.  Combined with (\ref{standingAssumptions})
we have  $n\le k$, which contradicts the hypothesis that $Y$ is a proper
submanifold of $X$.
\qed

\begin{corollary}\label{HardCor}
Let  $Y$ be a connected $k$-dimensional
submanifold of $X$, a connected projective manifold of
dimension $n\geq 2$.  Assume that $k>0$.  Assume that $L$ is a very ample 
line bundle on $X$
such that $(Y,L_Y)\cong\pnpair k 1$.   Let $\grs:\oline{X}\to X$ be the
blowing up of $X$ along $Y$.  Let $\phi:\oline{X}\to Z$
be the morphism associated to $L\otimes\sJ_Y$ as in {\rm
(\ref{SetUp})}. Then either $Z \cong \pn{n-k-1}$ with
$(X,L)\cong(\pn n,\sO_{\pn n}(1))$ or $\dim Z\geq k$.\end{corollary}
\proof It immediately follows by combining (\ref{Pk}) and
(\ref{ThmB}).\qed
\begin{corollary}\label{ceil}
Let  $Y$ be a connected $k$-dimensional
submanifold of $X$, a connected projective manifold of
dimension $n\geq 2$.  Assume that $k>0$.  Let $L$ be
a very ample line bundle on $X$. Assume that $L\otimes\sJ_Y$ is spanned by
its global sections.  Let $\grs:\oline{X}\to X$ be the blowing up of $X$
along $Y$.  Let $\phi:\oline{X}\to Z$ be the morphism
associated to $L\otimes\sJ_Y$ as in {\rm (\ref{SetUp})}. Assume that $Y$ is
not a complete intersection. Further assume that  $h^{2j}(Y,\rat)_{\rm
alg}=1$ for $j\leq n-k$. Then $\dim Z\geq n/2$.
\end{corollary}
\proof From (\ref{Dimension}) and (\ref{ThmB}) it follows that $\dim Z\geq
n-k$ and $\dim Z\geq k$. Thus $\dim Z\geq n/2$.\qed

Let us point out some relations between the results above and
Castelnuovo-Mumford regularity theory and the
Castelnuovo bound conjecture (see \cite{GLP}).

\begin{rem*}\label{CastelConj} Let $X$ be an $n$-dimensional smooth variety in
$\pn {n+r}$ and let $Y$ be a $k$-dimensional subvariety of $X$ of degree
$\delta:=L^k\cdot Y$ and $L:=\sO_{\pn {n+r}}(1)_X$.
  Assume that $k>0$. Let $q$ be the
codimension of $Y$ in  the smallest linear subspace $\pn {k+q}\subset\pn
{n+r}$ containing $Y$. The Castelnuovo bound conjecture says that
$(\delta-q+1)L\otimes\sJ_Y$ is spanned by global sections.
The conjecture is related to the question of
Castelnuovo-Mumford regularity, and it is known to hold when $Y$ has
dimension $1$ (see \cite{GLP}) or $2$ (see \cite{Laz}) and when $Y$ has
dimension $3$ and $\left[{{\delta-1}\over q}\right] \geq 6$ (see
\cite{Ran}). Assuming
the conjecture true and $Y$ a divisor, we will show that the
projection from $Y$ associated  to $L$ is birational except in certain
specific cases (see \S \ref{Divisorial}).\end{rem*}

\begin{rem*}\label{False} Let $X$ be an $n$-dimensional smooth variety in
$\pn {n+r}$ of degree $d$. Let $L$ be a very ample line bundle on $X$.
 From the Castelnuovo-Mumford regularity theory developed in \cite{Harris}
it follows that in case $n=1$,  for
$t>\left[\frac{d-1}{r}\right]$ one has $h^1(tL)=0$ and thus that
$h^0(tL)=\chi(\sO_X(tL))$.
 One might hope that this extends in
higher dimensions also. Unfortunately this is not true in dimension
$n\geq2$, as the following example shows.

Let $C_1$ be a smooth plane curve of degree $d_1$ with $L_1$ the
restriction of the hyperplane
section bundle of $\pn 2$ to $C_1$.
  Let $L_2$
be a very ample line bundle of degree $d_2:=d'+2g-2$, with $d'>0$, on  a
smooth curve $C_2$  of genus $g:=g(C_2)$.
Let $X:=C_1\times C_2\subset \pn N$ and let
$L:=p_1^*L_1\otimes p_2^*L_2$, where $p_i:X\to C_i$,
$i=1,2$, are the projections on the two factors. Note that if $g(C_i)\ge 2$ for
$i=1,2$, then  $X$ is a surface of general type. Since $d_2>2g-2$,
 we have $h^1(L_2)=0$ and
hence
$$h^0(L_2)=d_2-g+1=d'+g-1.$$
Thus by the  Kunn\"{e}th formula we have $h^0(L)=3(d'+g-1)$, i.e.,
$N=3(d'+g-1)-1$. Therefore $r=N-2=3(d'+g-2)$. Moreover
$d=2d_1d_2=2d_1(d'+2g-2)$. Thus the critical value, $c$, is
$c=\left[\frac{2d_1(d'+2g-2)-1}{3(d'+g-2)}\right]$. For a fixed $g$ and
taking
$d'\gg 0$ and $d_1 \ge 10$ we have
$\displaystyle c\sim\frac{2}{3}d_1<d_1-3.$
On the other hand, by using again Kunn\"{e}th formulas we get, for
$t=d_1-3$,\begin{equation}\label{nottrue}
h^1(tL)\geq h^1(\sO_{C_1}(t))=h^1(K_{C_1})=1.
\end{equation}
To show that the equality  $h^0(tL)=\chi(\sO_X(tL))$ for
$t>c=\left[\frac{d-1}{r}\right]$ is not true in general, consider the
smooth irreducible curve and the set $\Gamma$ of $d$ distinct points
obtained as transversal intersection of $X$ with a general $\pn {r+1}$ and
a general hyperplane $\pn r$ of the $\pn {r+1}$. Look at the exact sequence
$\ 0\to (t-1)L_C\to L_C\to tL_{\Gamma}\to 0,\ $
 From \cite[Theorem (3.7)]{Harris} we know that
$H^0(tL_C)\to H^0(tL_{\Gamma})$ is surjective for $t>c$ and therefore
$H^1((t-1)L_C)$ injects in $H^1(tL_C)$, for $t>c$. Since $h^1(tL_C)=0$ for
$t\gg 0$, we can conclude that
$h^1((t-1)L_C)=0$  for $t>c.$
Thus from the cohomology sequence associated to the exact sequence
$\ 0\to (t-2)L\to (t-1)L\to(t-1)L_C\to 0\ $ we infer that
$h^2((t-2)L)=0$ for $t>c.$
 From this we thus conclude that, for $t> c-2$, the equality
$h^0(tL)=\chi(\sO_X(tL))$  is equivalent to $h^1(tL)=0$. We have just
shown (see (\ref{nottrue})) that this is not the case.
\end{rem*}

\section{Examples}\label{ExSec}
\addtocounter{subsection}{1}\setcounter{theorem}{0}

In this section we give some examples to illustrate the results obtained in
\S \ref{StructureSec}. We use the same notation as in (\ref{SetUp}). The
following example shows that the dimension of the image of the projection
in Theorem (\ref{Dimension}) can actually reach all possible values.

\begin{example}\label{Bilbao2} Let $M$ be a smooth connected projective
variety of dimension $n-s$, $s\geq 0$. Let $X:=M\times\pn s$. Let $L$
be a very ample line bundle on $X$. Let $p:X\to \pn s$ be the product
projection and
set $H:=p^*\sO_{\pn s}(1)$. Let $Y$ be the $k$-dimensional subvariety of
$X$ obtained as transversal intersection of $n-k-1$ general members of
$|L|$ and a general $D_{n-k}\in |L-H|$. Assume that $k\geq s$ and $k>0$. Let
$\phi:\overline{X}\to Z$ be the morphism associated to $L\otimes\sJ$ as in
the usual set up (\ref{SetUp}).

Let
$N:=N_{Y/X}$ be the normal bundle of $Y$ in $X$. Note that
$N\cong(\oplus^{n-k-1}L_Y)\oplus(L-H)_Y$.
We let $V:=(\oplus^{n-k-1}\sO_X)\oplus H$ and $\sF:=(\oplus^{n-k-1}\sO_{\pn
s})\oplus\sO_{\pn s}(1)$. Thus $N^*(L)\cong V\cong p^*\sF$ and
  $E\cong{\Bbb
P}(p^*\sF)$, where $E$ is the exceptional divisor of the blowing up,
$\grs:\oline{X}\to X$, of $X$ along $Y$. Let $\gra$, $\grb$ be the
morphisms associated to $|\xi_{p^*\sF}|$
and $|\xi_{\sF}|$ respectively, where $\xi_{p^*\sF}$
and $\xi_{\sF}$ are the tautological line bundles of ${\Bbb P}(p^*\sF)$ and
${\Bbb P}(\sF)$.  Consider the projection
$p:X\to\pn s$. Since $\sF$ is a spanned vector bundle on $\pn s$, it is a
general fact that    $\gra:{\Bbb P}(p^*\sF)\to \pn
{N'}$ factors through   $\grb:{\Bbb P}(\sF)\to \pn
{N'}$.   Since $\sF$ is the direct sum of a
trivial bundle and a   very ample line
bundle, $\sO_{\pn s}(1)$,   $\xi_\sF$ is big. This
implies that
$\dim({\rm Im}\grb)=\dim({\Bbb P}(\sF))=n-k+s-1.$
Since $\dim\phi(\oline{X})\geq\dim\phi(E)=\dim Z$, it follows that
$\dim\phi(\oline{X})\geq n-k+s-1.$

Consider the Koszul complex
$$0\to\wedge^{n-k}V\otimes(-(n-k-1)L)\to\cdots\to\wedge^2V\otimes(-L)\to
V\to \sJ_Y(L)\to 0.$$
Set $T:=\oplus^{n-k-1}\sO_X$ and note that $\wedge^m(T\oplus H)=
\wedge^mT\oplus(\wedge^{m-1}T\otimes H)$ for each $m\geq 1$. Note also
that $h^0(H)=s+1$ and hence $h^0(V)=n-k+s$.
 From the hypercohomology sequence associated to the Koszul complex
above we see that
$h^0(\sJ_Y(L))=h^0(V)=n-k+s$. This is immediate if  $L-H$ is assumed ample, but
otherwise requires checking a few cases.  Thus  we conclude that the image
of the
morphism,
 $\phi:\oline{X}\to Z$, associated to $L\otimes\sJ_Y$ has dimension
$$
\dim\phi(\oline{X})=\dim Z\leq n-k+s-1\leq{\rm cod}_XY+s-1.
$$
Thus we conclude that $\dim\phi(\oline{X})=n-k+s-1$.

Note that the complete intersection situation  corresponds, in our
present notation, to the case $s=0$ with $p$ the constant map.
\end{example}

We have the following three infinite sequences of examples (for one more
class of examples see  (\ref{ExTh})  in \S \ref{Linear1}).

\begin{example}\label{Man7} (projection from a linear divisor) Let $X$ be
an $n$-dimensional projective submanifold of $\pn{2n-1}$.  Assume that there
is a linear $\pn{n-1}$, $D\subset X$.  Let $L$ denote the restriction of
$\pnsheaf {2n-1}1$ to $X$.  Since the morphism, $\psi :X\to \psi(X)$
associated 
to $|L-D|$ agrees with the restriction of the  projection of $\pn {2n-1}$
from $D$
away from $D$, we see that $\dim \psi(X)\le n-1$.  From this we conclude that
$(L-D)^n=0$.  A calculation given in Proposition (\ref{kn1})
shows that
$d:=L^n=\displaystyle\frac{(s+1)^n-1}{s}$ for $s\geq 1$ and $n$
for $s=0$, where the normal bundle of $D$ in $X$ is $\pnsheaf {n-1}{-s}$.
Since we have that $(L-D)_D\cong\pnsheaf {n-1}{s+1}$ is ample for $s\ge 0$,
we conclude
that if $s\ge 0$, then the morphism associated to
$|L-D|$ has at least an $(n-1)$-dimensional image.

We now show that such examples occur for all integers $n>0$ and $s\ge 0$.
Fix
integers $s\ge 0$ and $n>0$.
Let $\sP:=\proj{\pnsheaf {2n-1}1\oplus \pnsheaf {2n-1}{s+1}}$ and let $p
:\sP\to \pn {2n-1}$
denote the bundle projection.  Let $\xi$ denote the tautological line
bundle on $\sP$ such that
$p_*\xi\cong\pnsheaf {2n-1}1\oplus \pnsheaf {2n-1}{s+1}$. Note that by
counting constants we see
that the transversal intersection of $n$ general elements of
$|\xi |$ is a smooth $n$-fold $X'$ which maps isomorphically under $p$ to
its image $X$
in $\pn {2n-1}$.  Let $L:=\pnsheaf{2n-1}1_X$. Let $\sE:=\oplus^n\xi$. From
the Koszul complex
resolution of the ideal sheaf of $X'$ we get the exact sequence
$$0\to\det\sE^*\to\wedge^{n-1}\sE^*\to\cdots\to\wedge^2\sE^*\to\sE^*\to\sO_{\sP}
\to
\sO_{X'}\to 0.$$
By tensoring the sequence with $p^*\sO_{\pn {2n-1}}(1)$ we see that the
restriction map
gives an isomorphism
$H^0(\pn {2n-1},\sO_{\pn {2n-1}}(1))\cong H^0(X,L)$. Moreover the
intersection of $X'$ with the
section $\Sigma$ corresponding to the quotient
$\ \displaystyle \pnsheaf {2n-1}1\oplus \pnsheaf {2n-1}{s+1}\to\pnsheaf
{2n-1}1\ $
is a linear $\pn {n-1}$ with respect to $\pnsheaf {2n-1}1$. Thus $X'$
contains a linear $\pn {n-1}$.
Denote this by $D$. Since $N_{\Sigma/\sP}\cong\sO_{\pn
{2n-1}}(-s-1)\otimes\xi_\Sigma$ and $\xi_\Sigma\cong\pnsheaf {2n-1}1$, and
since the normal bundle $N_{D/X}$ of $D$ in $X$ is isomorphic to the
restriction of the normal bundle of $\Sigma$, we see that
$N_{D/X}\cong\sO_{\pn {n-1}}(-s)$.  As noted above the morphism,
$\phi:=p_{X}:X\to \pn {n-1}$,
 associated to $L\otimes\sJ_D$ has an $(n-1)$-dimensional image.

Recall that $L-D\approx \phi^*\sH$ for some ample and spanned line
bundle $\sH$ on $\pn {n-1}$. Then in the example above one has
$\ \displaystyle \sH^{n-1}=1$.
Indeed, let $\sH=\sO_{\pn {n-1}}(h)$. Since $L-D\approx \phi^*\sH$,
we see that $h^0(L-D)={{h+n-1}\choose{n-1}}$. From the exact sequence
$\ \displaystyle 0\to L-D\to L\to L_D\cong\sO_{\pn {n-1}}(1)\to 0\ $
we infer that
$h^0(L)\geq h^0(L-D)+n.$
Since $h^0(L)\leq 2n$ we conclude that $h^0(L-D)\leq n$. Thus, since
$n\geq 2$,
$\displaystyle {{h+n-1}\choose{n-1}}\leq n$ implies $h=1$.
\end{example}

The following example is related to Theorem (\ref{Alan}) in \S
\ref{Linear}.

\begin{example}\label{Hyper} We
construct here a smooth hypersurface of degree $d$ in $\pn {2k+1}$
containing a linear $\pn k$, such that the projection from the $\pn
k$ associated to $L:=\sO_X(1)$ has a $k$-dimensional image.

Consider in $\pn {2k+1}$ the degree
$d$ hypersurface defined by the equation
$$\sum_{j=0}^{2k+1}x_jx_{2k+1-j}^{d-1}=0.$$
Then $X$ is smooth and contains the linear $\pn k$ defined by the
equations $x_{2k+1}=\cdots=x_{k+1}=0$.  The projection from
this $\pn k$ has image $\pn k$.
\end{example}

\begin{example}\label{New} Let
$X:={\Bbb P}(\sE\oplus\sO_{\pn k}(1))$, where $\sE$ is a rank $r$
vector bundle on $\pn k$ of the form $\sE=\oplus_{i=1}^r\sO_{\pn
k}(a_i)$, $a_i\geq 1$. Then $X$ is of dimension $n=k+r$. Take as $\pn k$
the section of the $\pn r$-bundle $p:X\to \pn k$ corresponding to the
quotient
$$\sE\oplus\sO_{\pn k}(1)\to\sO_{\pn k}(1)\to 0.$$
This guarantees that $\xi_{\pn k}\approx \sO_{\pn k}(1)$, where
$\xi_{\pn k}$ is the restriction to $\pn k$ of the tautological
bundle $L:=\xi$ of $X$. Hence in particular $\delta:=L^k\cdot {\pn
k}=1$, i.e., $\pn k$ is linear.

Let $\grs:\oline{X}\to X$ be the blowing up of $X$ along the $\pn k$.
Note that $\grs$ induces the blowing up, $\pi:\oline{F}\to \pn r$, at
one point, $x$, of each fiber $F=\pn r$ of $p$. Consider the
morphism $\phi:\oline{X}\to Z$ associated  to $L\otimes\sJ_{\pn k}$.
Note that the restriction $\phi_{\oline{F}}$, for each fiber $F=\pn
r$, is the morphism given by the line bundle
$|\pi^*\sO_{\pn r}(1)-\pi^{-1}(x)|$. Therefore $\phi_{\oline{F}}$, being
the projection of $\pn r$ from the point $x$,  has lower
dimensional image. Since the fibers $F=\pn r$ cover $X$ we thus
conclude that $\phi$ has lower dimensional image.
\end{example}

\section{The divisorial case}\label{Divisorial}
\addtocounter{subsection}{1}\setcounter{theorem}{0}

In this section $L$  always denotes a very ample line bundle on a
$n$-dimensional projective manifold $X$, such that its global
sections, $\Gamma(L)$, embed $X$ in a projective space  $\pn {n+r}$.
Let $Y=D$ be a smooth connected divisor on $X$ of degree
$\delta=L^{n-1}\cdot D$.  We assume $n\geq 2$
since the case $n=1$ is trivial.

Recall that  $\delta L -  D$ is spanned (see Lemma
(\ref{Mumford})). In the present case we can
say considerably more. Let us first show the following fact.

\begin{lemma}\label{AmpleRestr} Let $L$ be a very ample line bundle
on a connected projective manifold $X$ of dimension $n$. Let $D$ be a
smooth divisor of degree $\delta=L^{n-1}\cdot D$. Then either
$(X,L,\sO_X(D))\cong (\pn n,\sO_{\pn n}(1),\sO_{\pn n}(\delta))$, or
the restriction $(\delta L-D)_D$ is an ample line bundle on
$D$.\end{lemma}
\proof By the conductor formula (\ref{ConductorF}) and the adjunction
formula we have that
\begin{equation}\label{Scerni}
(\delta -n-1)L_D-K_D\approx (\delta -n-1)L_D-(K_X+D)_D
\end{equation} is nef. By general adjunction theoretic results (see
e.g., \cite[(7.2.1)]{Book}) we know that $K_X+(n+1)L$ is either ample
or $(X,L)\cong(\pn n,\sO_{\pn n}(1))$. Therefore we see from
(\ref{Scerni}) that if $(\delta L-D)_D$ is not ample then
$K_{X|D}+(n+1)L_D$ is not ample and hence $(X,L)\cong (\pn n,\sO_{\pn
n}(1))$. In this case $\sO_X(D)\cong\sO_{\pn n}(\delta)$.\qed

Next, we recall the following definition.

\begin{def*}\label{kample} A line bundle, $L$, on a projective
variety, $X$, is $k$-ample for an integer $k\geq 0$, if $mL$ is
spanned for some $m>0$, and the morphism $X\to {\Bbb P}_{\comp}$
defined by $\Gamma(mL)$ for such an $m$ has all fibers of dimension
$\leq k$.\end{def*}

\begin{theorem}\label{oneAmple} Let $L$ be a very ample line bundle
on a connected projective manifold $X$ of dimension $n\geq 2$,
such that $\Gamma(L)$ embeds $X$ in $\pn {n+r}$. Let $D$ be a smooth
divisor on $X$ of degree $\delta=L^{n-1}\cdot D$. Then $\delta L-D$ is
$1$-ample except in the case when
$(X,L,\sO_X(D))\cong (\pn n,\sO_{\pn n}(1),\sO_{\pn n}(\delta))$.
\end{theorem}
\proof Let $F$ be a fiber of the morphism associated to
 $|\delta L-D|$ and assume $\dim F\geq 2$. Then  $(\delta
L-D)_F\approx \sO_F$, so that $D_F\approx\delta L_F$ is ample. This
implies that $D\cap F$ contains an effective curve, $C$, and $D\cdot
C>0$. But $(\delta L-D)\cdot C=0$ since $\delta L-D$ is trivial on
$F$. If  $(X,L,\sO_X(D))\not\cong (\pn n,\sO_{\pn n}(1),\sO_{\pn n}(\delta))$
this contradicts the ampleness of $(\delta L-D)_D$ (see
(\ref{AmpleRestr})).
\qed

If $\delta > 1$ we can say more.

\begin{theorem}\label{Properties} Let $L$ be a very ample line bundle
on a connected projective manifold $X$ of dimension $n\geq 2$,
such that $\Gamma(L)$ embeds $X$ in $\pn {n+r}$. Let $D$ be a smooth
divisor on $X$ of degree $\delta=L^{n-1}\cdot D>1$. Assume that
$(X,L,\sO_X(D))\not\cong (\pn n,\sO_{\pn n}(1),\sO_{\pn n}(\delta))$.
 Then the morphism associated to $|\delta L-D|$ is birational; moreover,
$\delta L-D$ is very ample if $n\geq r+2$.
\end{theorem}
\proof  First assume $n\geq r+2$, or, equivalently,   $2\dim X-(n+r)\geq
2$. Then
by the Barth-Lefschetz theorem (see e.g., \cite[(2.3.11)]{Book}) we
conclude that ${\rm Pic}(X)\cong \zed$ with generator the restriction of
the hyperplane section bundle on
 projective space. Since $\delta L-D$ is spanned
and not trivial unless $(X,L,\sO_X(D))\cong (\pn n,\sO_{\pn n}(1),\sO_{\pn
n}(\delta))$ (see (\ref{AmpleRestr})),
 we conclude $\delta L-D$ is a multiple of the restriction of the
hyperplane section bundle on projective space.
Thus $\delta L-D$ is very ample.

We next assume that $n \le r+1$. Then
$$\delta-1 > \frac{n}{r+1}(\delta-1)-n+1,$$
and Proposition (\ref{CastelP}) applies to say that $h^0((\delta -1)L-D) >
0$, from which it easily follows that
the morphism associated to $|\delta L-D|$ is birational.
\qed

Look at the embedding $X\subset\pn {n+r}$ and let $q$ be the codimension of
$D$ in
the smallest linear subspace $\pn {n-1+q}\subset\pn {n+r}$ containing
it. Let us assume that the Castelnuovo bound conjecture holds true, i.e.,
$(\delta -q+1)L-D$ is spanned by its global sections (compare with
(\ref{CastelConj})). Clearly we have
\begin{equation}\label{Uno}
r\geq q-1.
\end{equation}
Recall also the usual relations
\begin{equation}\label{Dos}
d\geq r+1
\ \ \mbox{\rm and\ \ }
\delta\geq q+1.
\end{equation}

\begin{prop}\label{PropertiesCast} Let $L$ be a very ample line bundle
on a connected projective manifold $X$ of dimension $n\geq 2$,
such that $\Gamma(L)$ embeds $X$ in $\pn {n+r}$. Let $D$ be a smooth
divisor on $X$ of degree $\delta=L^{n-1}\cdot D>1$. Let $q$ be the
codimension of $D$ in the smallest linear subspace
$\pn {n-1+q}\subset\pn {n+r}$ containing it. Assume that
$(\delta -q+1)L-D$ is spanned by its global sections. Then
\begin{enumerate}
\em\item\em If $n\geq r+2$, $(\delta-q+1)L-D$ is very ample unless
$X\cong\pn n$ and $\delta L\approx D$;
\em\item\em If $n\leq r+1$, then  the  morphism associated to
$|(\delta-q+1)L-D|$ is birational unless  $q=r+1$ and either
$n=r+1$  or $n<r+1$ and $\delta=r+2$.
\end{enumerate}\end{prop}
\proof Assume $n\geq r+2$. Let $d:=L^n$. We have the following fact.

\noindent {\bf Claim.\ }$(\delta-q+1)L-D$ is not trivial unless $X\cong \pn n$,
$\delta L\approx D$.

\noindent {\em Proof of  Claim.}
Assume $D\approx (\delta-q+1)L$. Dotting with $L^{n-1}$ gives
$(\delta-q+1)d=L^{n-1}\cdot D=\delta,$ or
$(d-1)\delta=d(q-1)$.
Using (\ref{Dos}) this gives
$(d-1)(q+1)\leq d(q-1),$ or
\begin{equation}\label{Four}
2d\leq q+1.
\end{equation}
Since by (\ref{Dos}) and (\ref{Uno}), $d\geq r+1\geq q$, we find $q\leq
1.$ Thus (\ref{Four}) yields $d=q=1$ and hence $r=0$ by (\ref{Dos}).
Therefore $X\cong \pn n$, $D\approx \delta L$.
$\Box$

Since $n\geq r+2$ is equivalent to $2\dim X-(n+r)\geq 2$, by the
Barth-Lefschetz theorem (see e.g., \cite[(2.3.11)]{Book}) we have ${\rm
Pic}(X)\cong \zed$. Since $(\delta-q+1)L-D$ is spanned and by the Claim
we can assume it is not trivial, we conclude that $(\delta-q+1)L-D$ is
very ample. This shows $1$).

As for $2$), assume $n\leq r+1$. If the morphism associated to
$|(\delta-q+1)L-D|$ is not birational, then $h^0((\delta-q)L-D)=0$.
Thus, by Proposition (\ref{CastelP}),
$\displaystyle \delta-q\leq\frac{n}{r+1}(\delta-1)-n+1,$ or
\begin{equation}\label{Five}
(\delta-1)\left(1-\frac{n}{r+1}\right)\leq q-n,
\end{equation}
or, by using $r\geq q-1$ from (\ref{Uno}),
\begin{equation}\label{Six}
(\delta-1)\left(\frac{r+1-n}{r+1}\right)\leq r+1-n.
\end{equation}
If $r+1=n$, then equality holds in (\ref{Six}) and hence in particular
$r=q-1$, i.e., $D$ spans $\pn {n+r}$. If $r+1>n$, inequality (\ref{Six})
yields $\delta-1\leq r+1$, or
$\delta\leq r+2.$
Since $q\leq\delta-1$ by (\ref{Dos}), inequality (\ref{Five}) gives
$$(\delta-1)\left(1-\frac{n}{r+1}\right)\leq \delta-1-n\ \ \mbox{\rm\  or\
}\ \
n\leq (\delta-1)\frac{n}{r+1}.$$
This implies $r+2\leq \delta$. Thus $\delta=r+2$. Also, at each step,
equalities hold true. Therefore $q=\delta-1=r+1$.\qed

\begin{example}\label{Range1} Notation as in (\ref{PropertiesCast}).
We give here an example in the range $n=r+1$ where 
 $|(\delta-q+1)L-D|$ is spanned but the morphism associated to it is not
birational, $D$ spans $\pn {n+r}$ and the projection from $D$ has an
$(n-1)$-dimensional image.

Consider the Segre embedding $X=\pn 1\times \pn {n-1}\hookrightarrow \pn
{n+r}=\pn {2n-1}$, $r=n-1$, and let $p_1:X\to \pn 1$, $p_2:X\to \pn
{n-1}$ be the projections on the two factors. Denote
$\sO(a,b):=p_1^*\sO_{\pn 1}(a)\otimes p_2^*\sO_{\pn {n-2}}(b)$, for given
integers $a$, $b$. Let $L:=\sO(1,1)$, so that $h^0(L)=2n$. Take a smooth
divisor $D$ in the linear system $|\sO(2,1)|$. We have
$\ d:=L^n=n\ $ and $\delta:=L^{n-1}\cdot D=n+1.$
Consider the exact sequence
$\ 0\to L-D\to L\to L_D\to 0.\ $
Note that $L-D=\sO(-1,0)$, so $h^0(L-D)=0$ and, by using Kunn\"{e}th's
formulas,
$h^1(L-D)=0$. Therefore $h^0(L)=h^0(L_D)$. This means that $D$ spans $\pn
{n+r}=\pn {2n-1}$, or $q=r+1=n$. Then
$(\delta-q+1)L-D=2L-D=\sO(0,1)$.
Thus $(\delta-q+1)L-D$ is not big, so that the projection from $D$
associated to it is not birational, and has an $(n-1)$-dimensional
image.
\end{example}

\begin{example}\label{Range2} Notation as in (\ref{PropertiesCast}).
We give here an example in the range $r=n$, where $(\delta-q+1)L-D$ is
spanned but not ample, in fact is $1$-ample, and the morphism associated
to it is birational.

Let $X:={\Bbb P}(\oplus^{n-1}\sO_{\pn 1}\oplus\sO_{\pn 1}(1))$. Let
$\xi$ be the tautological bundle of $X$ and let $F$ be a fiber of the
bundle projection $X\to \pn 1$. Let $L:=\xi+F$ and take a smooth divisor
$D\in |\xi+2F|$. Note that both $\xi+F$ and $\xi+2F$ are very ample (see
e.g., \cite[(3.2.4)]{Book}).

A standard check shows that $d=L^n=n+1$, $\delta=L^{n-1}\cdot D=n+2$ and
$h^0(L)=2n+1$, $h^0(L-D)=h^0(-F)=0$, $h^1(L-D)=0$. Thus $X\subset \pn
{2n}$, i.e., $q=r+1=n+1$. Then
$$(\delta-q+1)L-D=2L-D=\xi.$$
The line bundle $\xi$ is spanned but not ample (see e.g.,
\cite[(3.2.4)]{Book}) and the morphism associated to $|\xi|$ is the
blowing up $X\to \pn n$ of $\pn n$ along $\pn {n-2}$. Hence in particular
$\xi$ is $1$-ample.
\end{example}

\section{The linear case}\label{Linear}
\addtocounter{subsection}{1}\setcounter{theorem}{0}

Let $X$ be a smooth connected projective variety of dimension $n$,
polarized by a
very ample line bundle $L$. In this section we discuss some further results
about
the structure of projection maps from a
$k$-dimensional subvariety $Y$ of $X$, in the case when $Y$ is a linear $\pn k$
with respect to $L$. 

In (\ref{Alan}) we show that if the morphism, $\phi$, associated to
$L\otimes\sJ_Y$ as in  (\ref{SetUp}) has image dimension $n-k$, then
$\phi$ has $\pn {n-k}$ as image and $X$ is a hypersurface in $\pn {n+1}$.
Next we show in (\ref{HartC}) that
assuming ``Hartshorne's conjecture''  we have a stronger lower
bound for the dimension
of the image of $\phi$. Finally we prove in (\ref{Easy}) a
spannedness result for the adjoint bundle (see also (\ref{MoreAdj}) for more
adjunction theoretic structure type results in the case when $Y$ is a
codimension
$1$ linear $\pn {n-1})$.

Let us explicitly point out the following fact:
if  $Y$ is a smooth $k$-dimensional  subvariety of $(X,L)$ of degree
$\delta=L^k\cdot Y$, then, since
$\delta L\otimes\sJ_Y$  is spanned by global
sections by Lemma (\ref{Mumford}), the morphism associated to
$|tL\otimes\sJ_Y|$ is birational for $t\geq \delta+1$. In particular, if
$Y$ is a linear $\pn k$ with respect to $L$  and the projection
from $Y$ associated to $tL$  has lower dimensional image, then
necessarily $t=1$.

In the case when $Y$ is a linear $\pn k$ and the projection has image
dimension one bigger than the lowest possible value we have the following
result. We recall Theorem (\ref{Dimension}) for a general lower bound for
the image dimension of $\phi$ and we refer back to (\ref{Hyper}) which
gives in fact an example of the situation discussed below.

\begin{theorem}\label{Alan} Let $L$ be a very ample line bundle on $X$, a
connected projective manifold of dimension $n\geq 2$. Let $Y$ be a
subvariety of $X$ with $(Y,L_Y)\cong\pnpair k 1$.
 Let $\grs:\oline{X}\to X$ be the blowing up of $X$
along $Y$.  Let $\phi:\oline{X}\to Z$ be the morphism
associated to $L\otimes\sJ_Y$  as in
{\rm (\ref{SetUp})}. Assume that $\dim Z={\rm cod}_XY=n-k$ and $k\geq 2$.
Then $X$ is a hypersurface in $\pn{n+1}$.
\end{theorem}
\proof
Set $w:=n-k$. Since $\oline{L}-E$
is spanned and gives the projection  $\psi:\oline{X}\to
\psi(\overline{X})$ and since $\dim Z=w$, we have a surjection of locally
free sheaves $\oplus^{w+1}\sO_{\oline{X}}\to \oline{L}-E\to 0$. Hence,
restricting to $E$, we have an exact sequence
$$ 0\to \sK\to\oplus^{w+1}\sO_E\to
(\oline{L}-E)_E\to 0.$$ Consider the $\pn {w-1}$-bundle map $\pi:E\to
Y$. Let $N:=N_{X/Y}$ be the normal bundle of $Y$ in $X$.  Notice that
$(\oline{L}-E)_E\cong\xi$, the tautological line bundle of $E\cong{\Bbb
P}(N^*(L))$. By pushing forward under $\pi$, we get an exact sequence on
$Y$ \begin{equation}\label{ExSeq1}0\to K\to \oplus^{w+1}\sO_Y\to
N^*(L)\cong\pi_*\xi\to 0. \end{equation}
By comparing the ranks, since
$N^*(L)$ has rank ${\rm cod}_XY=w$, we conclude that $K$ is a line
bundle.

Since $Y\cong\pn k$, $k\geq 2$,   the first cohomology of a line
bundle is zero, i.e., $h^1(Y,K)=0$. This means that the
sections of $\oplus^{w+1}\sO_Y$ surject onto the sections of $N^*(L)$, so
$h^0(N^*(L))\leq w+1$.

Notice that $\oline{L}-E\approx\phi^*(\sH)$ for some ample line bundle
$\sH$ on $Z$. Since the restriction
$\phi_E:E\to Z$ is onto by (\ref{Structure}),  we have
$h^0(N^*(L))=h^0((\oline{L}-E)_E)\geq h^0(\sH)=h^0(\oline{L}-E).$
Thus
\begin{equation}\label{Mauro}
h^0(L\otimes\sJ_Y)=h^0(\oline{L}-E)\leq w+1.
\end{equation}
Now look at the exact sequence
$\ 0\to L\otimes \sJ_Y\to L\to L_Y\to 0.\ $
Recall that $L_Y\cong\sO_{\pn k}(1)$ since $Y$ is a linear $\pn k$.
Therefore, by (\ref{Mauro}),
$\displaystyle h^0(L)\leq h^0(L\otimes\sJ_Y)+h^0(\sO_{\pn k}(1))\leq
w+k+2=n+2$.
Thus, either $\Gamma(L)$ embeds $X$ as hypersurface in $\pn {n+1}$,
or else $h^0(L) = n+1$ and $X \cong \pn n$. However, the latter is
ruled out by the assumption $\dim Z \geq n-k$.
\qed

A minor modification of the proof of the theorem above gives us the
following result, which states that
assuming ``Hartshorne's conjecture'' (see \cite{HartshCI}) the image
dimension of $\phi$ has a stronger
lower bound unless $X$ is a complete intersection.

\begin{prop}\label{HartC} Let $X$ be a smooth connected projective variety
of dimension $n\geq 2$.
Let $L$ be a very ample line bundle on $X$. Let $Y$ be a linear $\pn k$
with respect to the embedding
given by $\Gamma(L)$. Let $\grs:\oline{X}\to X$ be the blowing up of $X$
along $Y$.
 Let $\phi:\oline{X}\to Z$ be the morphism associated to
$L\otimes\sJ_Y$  as in {\rm (\ref{SetUp})}.
 Assume  that Hartshorne's conjecture is true and that $X$ is not a
complete intersection. Then $\dim\phi(\oline{X})\geq {\rm
cod}_XY+\frac{k}{3}-1$.\end{prop}
\proof  First note that for $k\leq 2$ the bound on $\dim\phi(\oline{X})$
follows from Theorem
(\ref{Dimension}) and Corollary (\ref{Pk}), so we can assume $k\geq 3$.

Set $w:=n-k={\rm cod}_XY$ and $z:=\dim\phi(\oline{X})$. Exactly the same
argument as in the proof of Theorem
(\ref{Alan}) gives us an exact sequence 
$\ \displaystyle 0\to K\to \oplus^{z+1}\sO_Y\to N^*(tL)\to 0$ on $Y$,
where $K$ is a vector bundle of rank $z+1-w$ and $N$ is the normal bundle
of $Y$ in $X$.

Assume, by contradiction, that $z<w+\frac{k}{3}-1$, and therefore ${\rm
rank}(K)=z+1-w<\frac{k}{3}$. Thus from
(\ref{HartRev}) we know that $K$ splits as a direct sum of  line
bundles on $\pn k$ (here we are
using our present assumption that $k\geq 3$). Then the first cohomology of
$K$ is zero. This means that
the sections of $\oplus^{z+1}\sO_Y$ surject onto the sections of $N^*(L)$,
so $h^0(N^*(L))\leq z+1$.
Again, as in the proof of (\ref{Alan}), we thus conclude that
\begin{equation}\label{Mauro1}
h^0(L\otimes\sJ_Y)\leq z+1.
\end{equation}
Now look at the exact sequence
$\ 0\to L\otimes \sJ_Y\to L\to L_Y\to 0.\ $
Recall that $L_Y\cong\sO_{\pn k}(1)$ since $Y$ is a linear $\pn k$.
Therefore, by (\ref{Mauro1}),
$\ \displaystyle h^0(L)\leq h^0(L\otimes\sJ_Y)+h^0(\sO_{\pn k}(1))\leq z+k+2$.
Thus $\Gamma(L)$ embeds $X$  in $\pn {z+k+1}$. A direct numerical check
shows that the
inequality $z<w+\frac{k}{3}-1$ implies $n>\frac{2}{3}(z+k+1)$. Since we are
assuming that Hartshorne's conjecture is true, we thus conclude that $X$ is a
complete intersection.\qed

We need the following result.  The case when $k=1$ also follows
immediately from a result of Ilic \cite{Ilic}.
\begin{theorem}\label{newResult} Let $X$ be a connected projective
manifold of dimension $n\geq 2$.  Assume that $X$ is a $\pn {n-1}$-bundle
$\pi:X\to C$  over a smooth curve $C$ with fibers linear with respect
to $L$,  a very ample line bundle on $X$. Let $Y\subset X$ be a
linear $\pn k$  with respect to the embedding
given by $\Gamma(L)$.  Let $\grs:\oline{X}\to X$ be the blowing
up of $X$ along $Y$.
Let $\phi:\oline{X}\to Z$ be the morphism associated to $L\otimes\sJ_Y$
as in {\rm (\ref{SetUp})}. Then
$\dim Z<n$ if and only if either
\begin{enumerate}
\em\item\em $\dim \pi(Y)= 1$, $\dim Y=1$ and $Y$ is a section of $\pi$
corresponding to a surjection
from the vector bundle $\pi_*L$ onto a direct summand $\pnsheaf 1 1$; or
\em\item\em $\dim \pi(Y)=0$, $k=n-1$, and $(X,L)\cong (\pn {n-1}\times \pn
1,\sO_{\pn {n-1}\times \pn 1}(1,1))$.
\end{enumerate}
\end{theorem}
\proof We leave the reader to check the straightforward assertion that
$\dim Z<n$ in cases 1) and 2).
Assume now that $\dim Z<n$.

If $\dim \pi(Y)=1$, then since  $\pn k$ cannot map onto a curve if $k\ge
2$, we conclude that  $k=1$ and $C\cong \pn 1$.  Since $Y$ and
fibers of $\pi$ are linear, we conclude that $Y$ meets any given
fiber transversely in exactly one point.  Thus $Y$ corresponds to a
surjection $\pi_*L\to \pnsheaf 1 {L\cdot Y}\cong \pnsheaf 11$.
Using the fact that $\pi_*L$ is very ample and a direct sum of line bundles,
it is a simple check that $\pi_*L\to \pnsheaf 1 {L\cdot Y}\cong \pnsheaf 11$
splits.

Assume now that $\dim \pi(Y)=0$.  If the codimension of $Y$ is one, then
we have $0=(L-Y)^n=L^n-nL^{n-1}\cdot Y=L^n-n$.  From this we see that
$\pi_*L$ is a very ample rank $n$ vector bundle of degree $n$.  This
immediately implies that  $(X,L)\cong (\pn {n-1}\times \pn 1,\sO_{\pn
{n-1}\times \pn 1}(1,1))$.

Now we consider the case when the codimension of $Y$ is greater than one.
Since $N_{Y/X}\cong \sO_{\pn k}\oplus\oplus^{n-1-k}\pnsheaf k1$, it is a
straightforward
consequence of Lemma (\ref{differentials}) and   the fact that $N^*(L)$ is
spanned, that we
can choose $n-k-1$ smooth
divisors $D_1,\ldots,D_{n-k-1}$ in $|L\otimes\sJ_Y|$ all meeting
transversely in a smooth
$(k+1)$-dimensional subvariety $X_{k+1} := D_1 \cap \ldots \cap D_{n-k-1}$
containing $Y$
as a divisor.
But since it follows from the last paragraph that
$(X_{k+1},L_{X_{k+1}})\cong (\pn k\times\pn 1,\sO_{\pn k\times\pn
1}(1,1))$ we infer that $\pi_{X_{k+1}*}L_{X_{k+1}}\cong
\oplus^{k+1}\pnsheaf 1 1$. Thus  we conclude that $\ \displaystyle
\pi_*L\cong (\oplus^{n-k-1}\sO_{\pn 1})\oplus (\oplus^{k+1}\pnsheaf 1 1)$.
Since $\pi_*L$ is very ample, we conclude that $n=k+1$.
\qed

 The case $k=1$ of the following spannedness result for the 
adjoint bundle follows from \cite{Ilic}.

\begin{theorem}\label{Easy} Let $X$ be a smooth connected projective
variety of dimension $n\ge 2$. Let
$L$ be a very ample line bundle on $X$. Let $Y\subset X$ be a linear $\pn
k$ with respect to the embedding
given by $\Gamma(L)$.  Let $\grs:\oline{X}\to X$ be the blowing
up of $X$ along $Y$.
Let $\phi:\oline{X}\to Z$ be the morphism associated to $L\otimes\sJ_Y$
as in {\rm (\ref{SetUp})}.
Assume that $\dim Z<n$. Then
$K_X+(n-1)L$ is spanned by global sections unless either
\begin{enumerate}
 \em\item\em  $(X,L)\cong (\pn n,\sO_{\pn n}(1))$, $1\le k\le n-1$, with
$\dim Z=n-k-1$; or
\em \item\em  $(X,L)\cong(\sQ,\sO_{\sQ}(1))$, $\sQ$ a quadric in $\pn
{n+1}$, $1\le k\le \left[\frac{n}{2}\right]$, with $\dim Z=n-k$; or
\em \item\em  $(X,L)$ is a scroll, $\pi:X\to C$, over a smooth curve $C$,
i.e., $K_X+nL\approx \pi^*H$ for some ample
line bundle $H$ on $C$, with either
\begin{enumerate}
\item[{\rm (a)}] $\dim \pi(Y)= 1$, $\dim Y=1$ and $Y$ is a section of
$\pi$ corresponding a surjection
from the vector bundle $\pi_*L$ onto a direct summand $\pnsheaf 1 1$; or
\item[{\rm (b)}] $\dim \pi(Y)=0$, $k=n-1$, and $(X,L)\cong (\pn
{n-1}\times \pn 1,\sO_{\pn {n-1}\times \pn 1}(1,1))$.
\end{enumerate}
\end{enumerate}
\end{theorem}
\proof  From general adjunction theory results we know that
$K_X+(n-1)L$ is spanned unless either
\begin{enumerate}
\item[ (i)] $(X,L)\cong (\pn n,\sO_{\pn n}(1))$; or
\item[ (ii)] $(X,L)\cong(\sQ,\sO_{\sQ}(1))$, $\sQ$ a quadric in $\pn
{n+1}$; or \item[(iii)] $(X,L)$ is a scroll, $\pi:X\to C$, over a
smooth curve $C$, i.e., $K_X+nL\approx \pi^*H$ for some ample
line bundle $H$ on $C$.
\end{enumerate}
In case (i), by looking at the projection of $\pn n$ from $\pn k$, we see
that $2\le k\le n-1$ with $\dim
Z=n-k-1$. 

In case (ii) we see that $2\le k\le \left[\frac{n}{2}\right]$
with $\dim Z=n-k$, by looking at the
projection of $\pn {n+1}$ from $\pn k$.

In case (iii), use Theorem (\ref{newResult}).\qed

\begin{corollary}\label{Chern1} Let $X$ be a smooth connected projective
variety of dimension $n\geq 2$. Let
$L$ be a very ample line bundle on $X$. Let $Y\subset X$ be a linear $\pn
k$ with respect to the embedding
given by $\Gamma(L)$.  Let $\phi:\oline{X}\to Z$ be the
morphism associated to $L\otimes\sJ_Y$ as in {\rm
(\ref{SetUp})}. Assume that $\dim
Z<n$. Let $N:=N_{Y/X}$ be the normal bundle of $Y$ in $X$. If $(X,L)$ is
not as in one of cases $1)$, $2)$, $3)$ of {\rm (\ref{Easy})}, one has
$c_1(N)\leq n-2-k$. \end{corollary}
\proof By the assumption, $K_X+(n-1)L$ is spanned. On the other hand,
$$(K_X+(n-1)L)_Y\approx K_{Y}-\det N+(n-1)L_Y\cong\sO_{\pn
k}(n-2-k)-\det N.$$
Since $(K_X+(n-1)L)_Y\cong\sO_{\pn k}(b)$ for some nonnegative integer $b$,
we thus conclude that $\det N\cong
\sO_{\pn k}(a)$ for some integer  $a\leq n-2-k$.\qed

\section{The linear case in codimension $1$}\label{Linear1}
\addtocounter{subsection}{1}\setcounter{theorem}{0}

Let $X$ be a smooth connected projective variety of dimension $n\geq 2$.
Let $L$ be a very ample line
bundle on $X$. Let $P$ be a linear $\pn {n-1}\subset X$ with respect to
$L$, i.e., $\delta=L^{n-1}\cdot P=1$.
Recall that in this case the line bundle $L-P$ is spanned (see the
discussion after Lemma (\ref{imdiv})).
We follow the notation of (\ref{SetUp}), with the exception of denoting $Y$
by $P$ to
 emphasize its special nature.
 Thus we let $\psi :X\to \psi(X)$ be
the morphism associated to $|L-P|$ and $\psi=\frak{s}\circ \phi$ the Remmert-Stein
factorization of $\psi$ with
$\phi:X\to Z$ having connected fibers and $\frak{s}:Z\to\psi(X)$ finite.

In this section we study the projection  from $P$, a linear $\pn {n-1}$,
under the assumption that
$n>\dim\phi(X)$. For shortness, it is convenient to refer to the situation
above simply saying that
$(X,L,P)$ is a $\pn {n-1}$-{\em degenerate triple}.

First, let us state the following preliminary facts.

\begin{prop}\label{kn1} Let $X$ be a connected $n$-dimensional manifold and
let $L$ be very ample line
bundle on $X$. Assume that $(X,L,P)$ is a  $\pn {n-1}$-degenerate triple. Let
$N:=N_{\pn {n-1}/X}\cong\sO_{\pn
{n-1}}(-s)$ be the normal bundle of $P:=\pn {n-1}$ in $X$. Then we have:
\begin{enumerate}
\em\item\em $s\geq -1$, with equality only if $(X,L)\cong(\pn n,\sO_{\pn
n}(1))$;
\em\item\em  if $s\geq
0$, the morphism $\psi :X\to \psi(X)$ associated to $|L-P|$ has an
$(n-1)$-dimensional image with all fibers
having dimension one; and    $\psi_P$ is finite; and
\em\item\em the degree of $(X,L)$ is given by
$d:=L^n=\displaystyle\frac{(s+1)^n-1}{s}$ for $s\geq 1$ and by $n$
for $s=0$.
\end{enumerate}\end{prop}
\proof  Items 1) and 2) follow immediately from Lemma (\ref{AmpleRestr})
and Theorem (\ref{oneAmple}).

As for $3$), note that since $L-P$ is not big we have $(L-P)^n=0$. Then
$$d=\sum_{j=1}^n(-1)^{j-1}{{n}\choose{j}} L^{n-j}\cdot P^j.$$
By noting that
$L^{n-j}\cdot P^j=\sO_P(1)^{n-j}\cdot \sO_P(-s)^{j-1}=(-1)^{j-1}s^{j-1},$
we find
$\ \displaystyle d=\sum_{j=1}^n{{n}\choose{j}}s^{j-1}.$
This gives the result. \qed

In light of the above results we will assume that
$(X,L)\not\cong(\pn n,\sO_{\pn n}(1))$,
i.e., $N_{P/X}\cong\sO_{\pn {n-1}}(-s)$ with $s\geq 0$.

Let $\sH$ be
the ample line bundle on $Z$ such that $L-P\approx \phi^*(\sH)$. Set
$\frak{h}=\sH^{n-1}$ and $t=L\cdot f$
for a general fiber $f$ of $\phi$. We have
$L_P-P_P\approx\sO_{\pn {n-1}}(s+1)\approx \phi^*_P(\sH).$
Since $t=\deg\phi_P$, we conclude that
\begin{equation}\label{prodRel}
t\frak{h}=(s+1)^{n-1}.
\end{equation}
Note that the restriction $\phi_P:\pn {n-1}\to Z$ is a
$t$-to-one finite morphism.

\begin{rem*}Note that by
 (\ref{CM})  and  (\ref{Piclemma}) applied to the finite map $\phi_P$ we
 conclude that $Z$ is Cohen-Macaulay, has $t$-factorial
singularities,
and $\Pic Z\cong \zed$.
\end{rem*}

Let us give one more class of examples.

With the notation as
above, assume that $(X,L,P)$ is a $\pn {n-1}$-degenerate triple with  $s\geq
0$ and $t=1$. Since the restriction $\phi_P$  is an isomorphism
under this
assumption we see that,
by using also relation (\ref{prodRel}), $(Z,\sH)\cong (\pn {n-1},\sO_{\pn
{n-1}}(s+1))$, and that
$\phi$ is a $\pn 1$-bundle (see also \cite[(3.2.1)]{Book}). We let
$V:=\phi_*\sO_X(P)$ and
 \begin{equation}\label{VandE}
\sE:=\phi_*L\cong \phi_*(\sO_X(P)\otimes L)\cong\phi_*(\sO_X(P)\otimes
 \phi^*\sH) \cong V\otimes\pnsheaf{n-1}{s+1}.
\end{equation}
Then $X\cong\proj{\sE}\cong\proj{V}$.

\begin{prop}\label{ExTh} If $s\geq 0$ and $t=1$ then
$(X,L)\cong (\proj {\pnsheaf {n-1}{s+1}\oplus\pnsheaf{n-1}1},\xi)$, where
$\xi$
 denotes the tautological line bundle on $\proj {\pnsheaf
{n-1}{s+1}\oplus\pnsheaf{n-1}1}$.
\end{prop}
\proof From the exact sequence
$\ \displaystyle 0\to\sO_X\to\sO_X(P)\to P_P\cong\sO_{\pn {n-1}}(-s)\to 0,\ $
 by taking the direct image and since the higher direct image functor
$R^i\phi_*\sO_X$ is zero for $i>0$, we get the exact sequence
$0\to\sO_{\pn{n-1}}\to V\to \pnsheaf{n-1}{-s}\to 0.$
Since $h^1(\sO_{\pn {n-1}}(s))=0$ we see that this sequence splits. Thus
$\sE=\phi_*L\cong V\otimes\sO_{\pn {n-1}}(s+1)\cong{\pnsheaf
{n-1}{s+1}\oplus\pnsheaf{n-1}1}$. From
this the result is clear.\qed

 From relation (\ref{prodRel}) we see that $s=0$ implies $t=1$. This
gives the
following consequence.

\begin{corollary}\label{ExCor}If $s=0$ then
$(X,L)\cong (\pn{n-1}\times\pn1,\sO_{\pn{n-1}\times\pn1}(1,1))$.
\end{corollary}

\begin{rem*} Note that the example of a $\pn {n-1}$-degenerate
triple given by $\pn n$ blown
up at one point $z$, $p:X\to \pn n$, with $L=p^*\sO_{\pn n}(2)-P$,
$P=p^{-1}(z)$, fits in Proposition
(\ref{ExTh}) with $s=1$.\end{rem*}

By the above, we can   work from now on under the extra assumptions that
 $s\geq 1$ and $t\geq 2$,
where $N_{\pn {n-1}/X}\cong\sO_{\pn {n-1}}(-s)$ and $t=\deg\phi_{\pn
{n-1}}$.

We can now carry out some more adjunction theoretic analysis, improving, in
the case of a codimension $1$
linear projective space, the results proved in (\ref{Easy}). We will also
assume $n\geq 3$, since the problem
is completely solved when $n=2$ (see \cite{So1}, \cite[\S 8.4]{Book}). For
the structure of the first reduction
map occurring in the theorem below we refer to \cite[Chap. 7]{Book}.

\begin{theorem}\label{MoreAdj} Let $X$ be a smooth connected
$n$-dimensional variety, $n\geq 3$,  and let $L$
be very ample line bundle on $X$. Assume that $(X,L,P)$ is a  $\pn
{n-1}$-degenerate triple.
Let $N:=N_{\pn {n-1}/X}\cong\sO_{\pn {n-1}}(-s)$. Assume that $s\geq 1$ and
$t:=\deg \phi_P\geq 2$.
 Then the first
reduction exists, i.e., there exists a map $\pi:X\to X'$ expressing $X$ as
the blowup of a projective manifold $X'$ at a finite set $B$ with
$K_X+(n-1)L\approx \pi^* H$ for a very ample line bundle  $H$ on $X'$. Moreover it
follows that
$\pi$ is an isomorphism unless $B$ is a single point, $s=1$, and $P:=\pn
{n-1}=\pi^{-1}(B)$.
\end{theorem}
\proof  Set $P:=\pn {n-1}$. If $K_X+(n-1)L$ is not spanned,
then $(X,L)$ is as in one of cases $1$), $3$) of (\ref{Easy}) (notice that
case $2$) of (\ref{Easy}) is
excluded because we have $\dim P > \left[\frac{n}{2}\right]$). In
case $1$) we have that $K_X+(n+1)L$
is trivial, which implies that $\sO_P\approx(K_X+(n+1)L)_P\approx\sO_{\pn
{n-1}}(s+1)$. Thus $s=-1$. In case
$3$), we have $t=1$. Therefore
both cases $1$), $3$) of
(\ref{Easy}) are excluded in view of our present assumptions that $s\geq 1$
and $t\geq 2$.

Therefore we can assume that $K_X+(n-1)L$ is spanned.
It follows \cite[\S7.3]{Book}
that either $K_X+(n-1)L$ is nef and big  or:
\begin{enumerate}
\item $K_X\cong -(n-1)L$; or
\item $(X,L)$ is  a quadric fibration, $\pi:X\to C$,  over a  smooth curve
$C$, i.e., $K_X+(n-1)L\cong
\pi^*H$ for some ample line bundle $H$ on $C$; or
\item   $(X,L)$ is  a scroll, $\pi:X\to S$, over a smooth surface $S$,
i.e., $K_X+(n-1)L\cong
\pi^*H$ for some ample line bundle $H$ on $S$.
\end{enumerate}

In the first case we have that $\sO_f\cong (K_X+(n-1)L)_f$ for a general
fiber $f$ of $\phi$.
Since $(1-n)L\cdot f=(1-n)t=K_X\cdot f=\deg (K_f)$ we conclude that
$n=3$
and $t=L\cdot f=1$,
contradicting our present assumption $t\geq 2$.

Since $P=\pn {n-1}$ can't map to a curve by Lemma (\ref{BasicLemma}), we
conclude  in the second case
that $P$ is a component of a fiber of $\pi$. But since $n\geq 3$ fibers are
either irreducible
quadrics, or two $\pn {n-1}$'s meeting in a $\pn {n-2}$. Indeed multiple
fibers don't happen, since otherwise we could slice down to a surface and
have $\pn 1$ as a multiple fiber, which is a classical standard
impossibility. If we are in the case of two  $\pn {n-1}$'s meeting in a
$\pn {n-2}$, then we have negative normal bundle for each $\pn {n-1}$ and we
can contract one $\pn {n-1}$ to get a map of the other $\pn {n-1}$ to a
$(n-1)$-dimensional image but with the intersection $\pn {n-2}$ going to a
point, which is not possible again by Lemma (\ref{BasicLemma}).

In the third case we know from a result of the fourth author
\cite[Theorem (3.3)]{gottingen} that $\pi$ is
a $\pn {n-2}$-bundle.  Thus we conclude that $P$ is a section with $n=3$.
Indeed
since fibers of $\pi$ are one dimensional we conclude that $P$ meets a
general fiber $f$ of $\pi$ in a
finite nonempty set.  Since $L-P$ is nef and $L\cdot f=1$ we conclude that
$P\cdot f=1$.  Since
$(L-P)\cdot f=0$ it is clear that $\pi$ is the same as $\phi$ and
$t=1$.

Thus we see that $K_X+(n-1)L$ is big and the first reduction $\pi:X\to X'$
exists.
Assume that $\pi$ is not an isomorphism.  Let $F$ be a positive dimensional
fiber
of $\pi$.  We know that $F$ is a linear $\pn {n-1}$ with respect to $L$ and
$N_{F/X}\cong\sO_{\pn {n-1}}(-1)$.  If we show that
$F=P$ then we see that $s=1$ and the theorem will be proved.  Thus assume that
$F$ is not $P$.  Then we see that $F\cap P$ is empty or we would have the
absurdity that $\pi$
maps the positive dimensional subset $F\cap P$ of $P$ to the point $\pi(F)$
without mapping
$P$ to the same point.  Thus we have $ L_F\cong  \pnsheaf {n-1}
1$. Therefore
we see that $F$ is a section of $\phi: X\to Z$.  Thus we conclude that
$\phi$ is a
$\pn 1$-bundle over $\pn{n-1}$.  Restricting the bundle to a bundle $\phi_S
: S\to R$
on a smooth curve $R$ on $Z$, we find a $\pn 1$-bundle $S$ over $R$ with
two disjoint curves, $P \cap S$ and $F\cap S$, each with negative self
intersection since both the normal bundles $N_{P/X}$, $N_{F/X}$ are
negative. This is absurd.
\qed

We conclude this section by considering the special case of a threefold $X$.

\begin{prgrph*}{The three dimensional case}\label{X3} We use the same
notation  and assumptions as above. In
particular in view of the results above we make the blanket
assumption that  $s\geq 1$ and $t\geq 2$.

\begin{theorem}\label{threer1} Let $X$ be a smooth threefold and $L$ a very
ample line bundle
on $X$. Assume that $(X,L,P)$ is a $\pn 2$-degenerate triple. If $s=1$ and
${\frak t}:=\deg \phi_P\ge 2 $, then
$t=4$. In this case
$X$ is the blowing up at one point of the complete intersection of three
quadrics in $\pn
6$.\end{theorem}
\proof If $s=1$ then by Proposition (\ref{kn1}), $3$) we see that $L^3=7$.
Note
that we use the
classification of degree $7$ manifolds given in \cite{Io}. By Theorem
(\ref{MoreAdj}) we can assume
that $K_X+2L$ is nef and big. Thus quadric  fibrations over curves and
scrolls over curves and
surfaces are ruled out. By using the degree $7$ classification, two
possibilities remain.
\begin{enumerate}
\item $X$ is the blowing up at one point, $\pi:X\to X'$, of the complete
intersection $X'$ of three
quadrics in $\pn 6$, with $\pi$ the first reduction map; or
\item there exists a morphism $\rho:X\to C$ of $X$ to a curve $C$ given by
the complete linear system $|m(K_X+L)|$ for $m\gg 0$.
\end{enumerate}
In the first case we know from \cite{Io} that $L$ embeds $X$ into $\pn 5$.
This $X$ contains the
positive dimensional fiber of $\pi$ and thus since projection from this
linear $\pn 2$ must map to
$\pn 2$ we conclude that this is an example with $s=1$. Let $f\cong\pn 1$
be a fiber of $\phi:X\to
Z$. To see what $t$ is, note that
$K_X+L$ being  nef yields $t=L\cdot f\geq -K_X\cdot f=2$. By Theorem
(\ref{MoreAdj}) we know that
$P$ coincides with the exceptional divisor of $\pi$. Moreover,
$-K_{X'}\cong\sO_{X'}(1)=L'$, the
polarization of the first reduction $X'$, which satisfies the condition
$L\cong\pi^*L'-P$. Then
$$K_X\cong\pi^*K_{X'}+2P\cong -L-P+2P=-L.$$
Hence we have $K_X\cdot f=\deg(K_f)=0$. Thus we cannot have $t=2$ since
this would imply $f$ was
rational. Since we are assuming $t\geq 2$ we conclude by relation
(\ref{prodRel}) that $t$ must equal
$4$.

In the second case $\rho(P)$ must be a point by Lemma (\ref{BasicLemma})
and therefore
$(K_X+L)_P\cong\sO_P$. Since $(K_X+L)_P\cong\sO_P(s-2)$ we get the
contradiction $s=2$.\qed

Combining Theorem (\ref{MoreAdj}) and Theorem (\ref{threer1}) we have the
following result.

\begin{corollary}\label{firstRed1}Let $X$ be a smooth threefold and $L$ a
very ample line
bundle on $X$. Assume that $(X,L,P)$ is a $\pn 2$-degenerate triple. If
$s\ge1$ and
${\frak t}:=\deg \phi_P\ge 2 $, then
either $X$ is the blowing up at one point of the
complete intersection of three quadrics in $\pn 6$, or $K_X+2L$ is very
ample.\end{corollary}
\proof By (\ref{MoreAdj}) and  (\ref{threer1}) we know that either $s=1$
and $X$ is the blowing up at one point of the  complete
intersection of three quadrics in $\pn 6$ or $X$ is isomorphic to its own
first reduction.\qed

\begin{theorem}\label{caset2}Let $X$ be a smooth threefold and $L$ a very
ample line bundle
on $X$. Assume that $(X,L,P)$ is a $\pn 2$-degenerate triple. Further assume
$s\geq 2$. Then the
case $t=2$ does not occur.\end{theorem}
\proof By Corollary (\ref{firstRed1}) we can assume that $(X,L)$ is its own
first reduction. A simple
check of the list of pairs with $K_X+L$ not nef (see \cite[\S7.3]{Book})
shows that they cannot
occur if $s\geq 2$. Thus we can assume that $K_X+L$ is nef. We know that
there is a morphism
with connected fibers $\rho:X\to W$ of $X$ onto a normal variety $W$, given
by $|m(K_X+L)|$ for
$m\gg 0$, with $K_X+L\cong\rho^*H$ for some ample line bundle $H$ on $W$.
Note that if $t=2$ then
the general fiber of $\phi:X\to Z$ is a conic. Thus $K_X+L$ must be trivial
on the general fiber of
$\phi$. Then there exists a surjective morphism $q:Z\to W$ such that
$q\circ\phi=\rho$, whence $\dim
W\leq 2$. Note also that $\dim W>0$. Indeed otherwise $K_X+L$ would be
trivial and therefore, since
$(K_X+L)_P\cong\sO_P(s-2)$, we would have $s=2$. But $t=s=2$
contradicts
relation (\ref{prodRel}).

The divisor $P$ can not be in a fiber of $\rho$. If it was we would have
$(K_X+L)_P\cong\sO_P$. This
would imply $s=2$. Then again $t=s=2$ contradicts relation
(\ref{prodRel}).
By using Lemma
(\ref{BasicLemma}) we conclude that $\dim W=2$ and, since $P$ must map onto
$W$, that all fibers of
$\rho$ are one dimensional. By the above, $(X,L)$ is a quadric fibration over
the surface $W$. Then by
Besana's results \cite{besana} we know that $W$ is smooth and thus by
Lazarsfeld's theorem (see
e.g., \cite[(3.1.7)]{Book}) we know that $W$ is $\pn 2$. We also see that
the maps $\rho$ and $\phi$ are the same.

Note that by pulling back to $P$ we have
$$m(K_X+L)_P\cong\sO_P(m(s-2))\cong(L-P)_P\cong\phi_P^*\sH\cong\sO_{\pn
2}(s+1).$$
This gives $s+1=m(s-2)$ and hence either $s=5$, $m=2$, $L-P\cong
2(K_X+L)$, or $s=3$, $m=4$, $L-P\cong 4(K_X+L)$. Assume $s=5$. Then, since
$t=2$, relation (\ref{prodRel}) gives $\frak h=\sH^2=18$. But since
$L-P\cong2(K_X+L)$ we have the absurdity that $18=\sH^2=4H^2$. Assume $s=3$.
Then $\frak h=\sH^2=8$ from relation (\ref{prodRel}) and $L-P\cong 4(K_X+L)$
gives the absurdity $8=\sH^2=16H^2$.\qed
\end{prgrph*}

\bigskip
{
\begin{tabular}{l} Mauro C. Beltrametti\\
 Dipartimento di Matematica\\
Via Dodecaneso 35\\
 I-16146 Genova, Italy\\
beltrame@dima.unige.it\\

\ \ \ \\
Alan Howard\\
Department of Mathematics\\
Notre Dame, Indiana, 46556, U.S.A.\\
howard.1@nd.edu\\

\ \ \ \\
Andrew J. Sommese\\
Department of Mathematics\\
Notre Dame, Indiana, 46556, U.S.A.\\
sommese.1@nd.edu\\
\verb+URL:\ \ http://www.nd.edu/~sommese/index.html +\\

 \end{tabular}
}

\end{document}